\documentclass[12pt]{article}
\usepackage{amssymb,amsmath,bm}
\usepackage{mathrsfs}
\usepackage{natbib}
\usepackage{constants}
\begin{document}
\newcommand{\bea}{\begin{eqnarray}}
\newcommand{\ena}{\end{eqnarray}}
\newcommand{\beas}{\begin{eqnarray*}}
\newcommand{\enas}{\end{eqnarray*}}
\newcommand{\beq}{\begin{equation}}
\newcommand{\enq}{\end{equation}}
\def\qed{\hfill \mbox{\rule{0.5em}{0.5em}}}
\newcommand{\bbox}{\hfill $\Box$}
\newcommand{\ignore}[1]{}
\newcommand{\ignorex}[1]{#1}
\newcommand{\abs}[1]{\left\vert#1\right\vert}
\newcommand{\wtilde}[1]{\widetilde{#1}}
\newcommand{\qmq}[1]{\quad\mbox{#1}\quad}
\newcommand{\qm}[1]{\quad\mbox{#1}}
\newcommand{\var}{\mbox{Var}}
\newcommand{\nn}{\nonumber}
\newcommand{\Bvert}{\left\vert\vphantom{\frac{1}{1}}\right.}
\newcommand{\To}{\rightarrow}
\newcommand{\eps}{\varepsilon}

\newtheorem{theorem}{Theorem}[section]
\newtheorem{corollary}{Corollary}[section]
\newtheorem{conjecture}{Conjecture}[section]
\newtheorem{proposition}{Proposition}[section]
\newtheorem{lemma}{Lemma}[section]
\newtheorem{definition}{Definition}[section]
\newtheorem{example}{Example}[section]
\newtheorem{remark}{Remark}[section]
\newtheorem{case}{Case}[section]
\newtheorem{condition}{Condition}[section]
\newcommand{\pf}{\noindent {\it Proof:} }
\newcommand{\proof}{\noindent {\it Proof:} }

\title{{\bf\Large A Berry-Esseen bound for the uniform multinomial occupancy model}}
\author{Jay Bartroff and Larry Goldstein\thanks{Both authors work partially supported by NSA
grant H98230-11-1-0162.}\\University of Southern California}
\footnotetext{MSC 2010 subject classifications: Primary 60F05\ignore{Central limit and other weak theorems},
60C05\ignore{Combinatorial probability}}
\footnotetext{Key words and phrases: Stein's method, size bias, coupling, urn models}
\maketitle
\date{}

\begin{abstract}
The inductive size bias coupling technique and Stein's method yield a Berry-Esseen theorem for the number of urns having occupancy $d \ge 2$ when $n$ balls are uniformly distributed over $m$ urns. In particular,
there exists a constant $C$ depending only on $d$ such that
$$\sup_{z \in \mathbb{R}}\left|P\left( W_{n,m} \le z \right) -P(Z \le z)\right| \le C \left( \frac{1+(\frac{n}{m})^3}{\sigma_{n,m}}\right) \quad\mbox{for all $n \ge d$ and $m \ge 2$,}$$
where $W_{n,m}$ and $\sigma_{n,m}^2$ are the standardized count and variance, respectively, of the number of urns with $d$ balls, and $Z$ is a standard normal random variable. Asymptotically, the bound is optimal up to constants if $n$ and $m$ tend to infinity together in a way such that $n/m$ stays bounded.
\end{abstract}

\section{Introduction}
\label{sec:intro}

In this paper we provide a Berry-Esseen theorem in the classical occupancy problem for the normal approximation of the distribution of the number of urns having occupancy $d$ when $n$ balls are uniformly distributed among $m$ urns. Our proof relies on the inductive version of Stein's method using size bias couplings as presented in \citet{Go12}. In turn, that work springs from the use of induction in \citet{Bolthausen84}, achieving bounds for the combinatorial central limit theorem. The inductive method relies on expressing a bound for the distance of the given variable to the normal in terms of smaller versions of the same problem. For instance, in the occupancy model, conditional on the contents of a randomly chosen urn, the distribution of the remaining balls is uniformly distributed over one fewer urn.

Stein's method often proceeds by coupling a random variable $Y$ of interest to a related one using, for example, the method of exchangeable pairs, size bias couplings, or zero bias couplings (see \citet{Stein72}, \citet{Stein86} and \citet{Chen10}). However, some of the couplings that are the simplest to construct may lack a key boundedness property that is required for the application of many results. By applying a theorem that does not require the coupling to be bounded, in Theorem \ref{thm:occ} we are able to extend the work of \citet{Englund81} on the number of empty urns, and that of \citet{Penrose09} on the number of urns occupied by a single ball, to the case of all occupancies of size two and greater.

In the general multinomial occupancy model, one considers a vector ${\bf M}_n$ having components $M_n(i)$ that record the number of balls falling in urn $i$ in $n$ independent trials, where in each trial a single ball falls in urn $i$ with probability $\theta_i$ for all $i \ge 1$. In particular, the (multinomial) distribution ${\cal M}(n,\theta)$ of ${\bf M}_n$ is given by
\beas 
P(M_n(i)=m_i, i \ge 1) = \frac{n!}{\prod_{i \ge 1}m_i!} \prod_{i \ge 1} \theta_i^{m_i}
\enas
when $m_i$, $i \ge 1$, are nonnegative integers summing to $n$, and $\theta \in \Theta$ where
\begin{equation}
\label{def:Theta}
\Theta=\left\{(\theta_1,\theta_2,\ldots): \theta_i \ge 0, i \ge 1, \sum_{i \ge 1} \theta_i=1\right\}.
\end{equation}
For all $d \ge 0$ the number $Y_n^{(d)}$ of urns containing $d$ balls is therefore given by
\begin{equation}
 \label{Yn:balls}
Y_n^{(d)}= \sum_{i \ge 1} X_{n,i}^{(d)} \qmq{where} X_{n,i}^{(d)}={\bf 1}(M_n(i)=d).
\end{equation}

Among the many applications of multinomial occupancy models are the well-known species trapping problem (see \citet{Chao96}, \citet{Robbins68b}, or \citet{Starr79}) and the closely-related problem of statistical linguistics (see \citet{Efron76} and \citet{Thisted87}). In these applications a collection of species are trapped, or a collection of words are observed, according to the multinomial distribution ${\cal M}(n,\theta)$, and estimators of parameters related to the number of unseen species, or words known but unused by an author, are of central interest. Estimators of, say, the number of unknown species most often take the form of linear combinations of $Y_{n}^{(d)}$ for various $d$. For example, a well-known conjecture of \citet{Starr79} is that the uniformly minimum variance unbiased estimator, or UMVUE, of the probability of sampling a new species in a sample of size $n-n_0$, based on an original sample of size $n_0$, is
\beas 
\sum_{d=1}^{n-n_0} \frac{{n-n_0-1\choose d-1}}{{n\choose d}}Y_{n}^{(d)}.
\enas

For occupancy models where $n$ balls are distributed among the first $m$ urns only, the urn probability vector $\theta$ is given by $(\theta_1, \ldots, \theta_m,0,0,\ldots) \in \Theta$ as in (\ref{def:Theta}). Below we will find it convenient to continue to consider the case where the urns are indexed by all $i \ge 1$, even though all but the first $m$ of them will be empty.

In what follows we fix $d\ge 0$ and drop the superscript $(d)$ from our notation, denoting $X_{n,i}^{(d)}$ and $Y_n^{(d)}$ simply as $X_{n,i}$ and $Y_n$, respectively. \citet[][p.~37]{Kolchin78} show that the mean $\mu_{n,m}$ and variance $\sigma_{n,m}^2$ of the number $Y_n$ of urns occupied by $d \in \{0,1,\ldots\}$ balls, when $n$ balls are distributed uniformly over $m$ urns, are given by
\begin{align}
\mu_{n,m}&=m{n \choose d}\frac{1}{m^d}\left(1-\frac{1}{m}\right)^{n-d}, \quad \mbox{and} \label{occ:mv}\\
\sigma_{n,m}^2 &=\mu_{n,m} + m(m-1){n \choose d,d,n-2d}\frac{1}{m^{2d}}\left( 1 - \frac{2}{m} \right)^{n-2d}-\mu_{n,m}^2\label{varYn}
\end{align}
for all $n \ge d$ and $m \ge 2$, with the second term in (\ref{varYn}) set to zero for all $d \le n < 2d$ and $m \ge 2$.

Since the cases $d=0$ and $d=1$ having already been handled by \citet{Englund81} and \citet{Penrose09}, respectively, we focus on $d \ge 2$. Our main result is the following.

\begin{theorem} \label{thm:occ}
For $d \in \{2,3,\ldots,\}$, let $Y_n$ be the number of urns containing $d$ balls in the uniform occupancy model with $n$ balls and $m$ urns. Then, with $\mu_{n,m}$ given by (\ref{occ:mv}), $\sigma_{n,m}^2$ by (\ref{varYn}), and
\bea \label{rnm:def}
r_{n,m}=\frac{\sigma_{n,m}}{1+ (n/m)^3},
\ena
there exists a constant $C$ depending only on $d$ such that the standardized count
\beas
W_{n,m}=\frac{Y_n-\mu_{n,m}}{\sigma_{n,m}}
\enas
satisfies
\begin{equation}\label{main.thm:sup.bd}
\sup_{z \in \mathbb{R}}|P(W_{n,m} \le z)-P(Z \le z)| \le C/r_{n,m} \quad \mbox{for all $n \ge d$ and $m \ge 2$.}
\end{equation}
\end{theorem}

Regarding lower bounds, \citet[][Section~6]{Englund81} shows that in the case $d=0$,
\begin{equation}
\label{englund.lower}
\sup_{z \in \mathbb{R}}|P(W_{n,m} \le z)-P(Z \le z)| \ge 0.087/\max(3,\sigma_{n,m})
\end{equation}
and we remark that Englund's argument holds without changes for any random variable $W_{n,m}$ with finite variance supported on the integers, and  so for the $d \ge 0$ cases of the occupancy problem in particular.

Although Theorem \ref{thm:occ} yields a bound for all $n\ge d$ and $m\ge 2$, often interest centers on the behavior of a sequence of occupancy models where $n$ and $m$ vary together in such a way that the ratio of $n$ to $m$ is bounded away from zero and infinity, that is, when there exist $0<a<b<\infty$ such that
\bea \label{central.domain.a.b}
a \le \frac{n}{m} \le b.
\ena
Note that the bound on the supremum norm in \eqref{main.thm:sup.bd} achieves the rate $1/\sigma_{n,m}$, optimal in view of (\ref{englund.lower}), when $\sigma_{n,m}/r_{n,m}$ is bounded away from infinity, or equivalently, when the upper bound in (\ref{central.domain.a.b}) holds.
These observations yield the following immediate corollary; see also Section~\ref{discussion} for a more detailed discussion of these and other asymptotic regimes.
\begin{corollary}\label{cor:occ}
For any $b \in (0,\infty)$ there exists a constant $C$, depending only on $d \in \{2,3,\ldots\}$ and $b$, such that
\beas
\sup_{z \in \mathbb{R}}|P(W_{n,m} \le z)-P(Z \le z)| \le C/\sigma_{n,m}
\enas
for all $n \ge d$ and $m \ge 2$ that satisfy $n/m \le b$, and the bound is optimal up to constants.
\end{corollary}
Specializing the broad results of Theorem 4.2 of \citet{ChRo10} for general functions of urn occupancies to the case considered here under (\ref{central.domain.a.b}) results in a bound in Kolmogorov distance in the central domain such as the one here, with explicit constants but additional factors of $\log n$ to various powers.

To begin to describe the first ingredient required for the proof of Theorem \ref{thm:occ}, the construction of a size biased coupling, recall that for a nonnegative random variable $Y$ with finite, nonzero mean $\mu$, $Y^s$ has the $Y$-size bias distribution if
\bea \label{intro:szb}
E[Yf(Y)]=\mu Ef(Y^s)
\ena
for functions $f$ for which the expectations above exist. In employing the size bias version of Stein's method, see \citet{Baldi89}, \citet{Goldstein96} and \cite{Chen10}, the goal is to construct, on the same space as $Y$, a variable $Y^s$ with the $Y$-size bias distribution, such that $Y$ and $Y^s$ are close is some sense. Previous applications of size bias coupling in Stein's method for producing bounds in the Kolmogorov distance in the presence of dependence, but for \citet{Go12}, have required that $|Y^s-Y|$ be bounded.

To size bias the number of urns $Y_n$ containing $d$ balls, note that
when $n$ balls are uniformly distributed over $m$ urns, $Y_n$ in (\ref{Yn:balls})
is the sum of $m$ exchangeable indicators. In general, Lemma \ref{sblem} below says, essentially, that to size bias such a sum one chooses an indicator uniformly, sets it to one, and then `adjusts' the remaining indicators, if necessary, to have their original distribution given that the selected indicator now takes the value one. In the occupancy problem, to set an indicator for a chosen urn equal to one if it is not so already, one must either add balls to that urn if it has fewer than $d$ balls, or redistribute balls from the urn if it has an excess over $d$. As it is possible that the chosen urn has, say, all $n$ balls, the resulting coupling fails to be bounded in $n$. However, as there is small probability that a very large number of balls will need to be redistributed, the coupling can be controlled using quantities such as moments on bounds $K_n$ on the absolute difference between $Y_n^s$ and $Y_n$.

To describe the second ingredient in the application of Theorem \ref{thm:main} in general, the inductive component, suppose that for some nonnegative integer $n_1$, for all $n \ge n_1$ we are given a nonnegative random variable $Y_n$ whose distribution ${\cal L}_\theta$ depends on a parameter $\theta$ in a topological space $\Theta_n$. As bounds to the normal for $Y_n$ can be expressed in terms of a number of quantities, including bounds to the normal for `smaller versions' of the same problem, an inductive argument yielding a recursion for the bound may be constructed 
when for random variables $L_n$ and $\psi_{n,\theta}$ taking values in
 $\{0,\ldots,n\}$ and $\Theta_{n-L_n}$, respectively, and a certain collection of random variables $J_n$ there exists a random variable $V_n$ on the same space as $Y_n$ such that
\beas 
{\cal L}_\theta(V_n|J_n)={\cal L}_{\psi_{n,\theta}}(Y_{n-L_n})
\enas
holds on a set where the size of $L_n$ is controlled.
One must also control the difference between $Y_n$ and $V_n$, but again strict boundedness is not required on $Y_n-V_n$ but rather moment estimates of a bounding random variable $B_n$ satisfying $|Y_n-V_n| \le B_n$.

Regarding the inductive component for our occupancy problem, if the urn chosen to have occupancy $d$ in the size bias configuration is removed, then, conditional on the identity of that urn and the number of balls it contains, the remaining configuration has the same uniform multinomial distribution over the remaining urns, one fewer than the number in the original configuration, of the balls not contained in the urn chosen. And again, as with the bound $K_n$ on $|Y_n^s-Y_n|$, though it is possible that the chosen urn contains a very large number of balls, it is unlikely that it will.

In the uniform model, \citet{Englund81} gave an explicit Berry-Esseen bound of order $1/\sigma_{n,m}$, with a corresponding lower bound (\ref{englund.lower}) of the same order, for the number of occupied urns, or equivalently, for the number of empty urns forming the complement, that is, those with occupancy $d=0$. For the non-uniform case, \citet{Quine84} gave a less explicit error bound. \citet{Hwang08} obtained a local limit theorem, and also describe applications including species trapping and statistical linguistics.
\citet{Johnson77} and \citet{Kolchin78} give results for models of this type in the uniform and some non-uniform cases. \citet{Penrose09} considers the case $d=1$ where $Y_n$ counts the number of isolated balls, and obtains a Berry-Esseen bound via size-biased coupling in the uniform case, and for the non-uniform case as well with a slightly larger constant. \citet{Karlin67}, \citet{Gnedin07} and \citet{Barbour09} consider the infinite occupancy model, the first two proving central limit theorems for the number of occupied urns, the last providing a multivariate normal approximation for arbitrary occupancies of a fixed number of urns.

In Section \ref{size.bias.coupling} we construct the coupling of $Y_n$ and the size biased variable $Y_n^s$. In Section \ref{proof.of.theorem}, with the help of Lemma \ref{r1.large}, we prove Theorem \ref{thm:occ} by verifying the conditions of Theorem \ref{thm:main}. Some discussion is provided in Section \ref{discussion}, and the proof of Lemma \ref{r1.large} is given in the Appendix. With $\mathbb{Z}$ the set of integers, let $\mathbb{N}_k=\mathbb{Z} \cap [k,\infty)$. Throughout, we will use $C,C_1,C_2,\ldots$ to denote positive, finite constants depending only on $d$. Since in what follows we focus on the uniform occupancy problem, for notational simplicity we specify the multinomial probability vector by $m\in \mathbb{N}_1$ rather than by the corresponding vector
\begin{equation}\label{theta_m}
\theta_m=(\underbrace{1/m,\ldots,1/m,}_{m}0,0,\ldots)
\end{equation}
and write $\mathbb{N}_2$ for our parameter space.  When considering subsets $\Theta_n \subseteq \Theta$ for some $n \in \mathbb{N}$ and invoking Theorem \ref{thm:main}, statements such as $m \in \Theta_n$ should be interpreted as meaning that $\theta_m \in \Theta_n$. Further, we will denote the uniform multinomial distribution of $n$ balls over $m$ urns as ${\cal M}(n,m)$, in parallel to our notation for the binomial ${\cal B}(n,p)$ distribution with $n$ trials and success probability $p$. For ${\bf M}_n \sim {\cal M}(n,m)$, in accordance with (\ref{theta_m}), we have $M_n(j)=0$ for all $j > m$. 

\section{Size Bias Coupling}
\label{size.bias.coupling}
A general prescription for size biasing a sum of nonnegative variables is given in \citet{Goldstein96}; specializing to exchangeable indicators yields the following result.
\begin{lemma} \label{sblem}
Suppose $Y=\sum_{\alpha \in {\cal I}}X_\alpha$, a finite sum of nontrivial exchangeable
Bernoulli variables $\{X_\alpha, \alpha \in {\cal I}\}$, and that for $\alpha \in {\cal I}$
the variables
$\{X_\beta^\alpha, \beta \in {\cal I}\}$ have joint distribution
\beas 
{\cal L}(X_\beta^\alpha, \beta \in {\cal I})={\cal L}(X_\beta, \beta \in {\cal I}|X_\alpha=1).
\enas
Then
$$
Y^\alpha=\sum_{\beta \in {\cal I}} X_\beta^\alpha
$$
has the $Y$ size biased distribution $Y^s$ characterized by (\ref{intro:szb}), as does the mixture $Y^I$ when $I$ is a random index with values in ${\cal I}$, independent of all other variables.
\end{lemma}

\pf
First, fixing $\alpha \in {\cal I}$, we show that $Y^\alpha$ satisfies (\ref{intro:szb}). For given $f$,
\beas
E[Y f(Y)] = \sum_{\beta \in {\cal I}} E[X_\beta f(Y)]
 = \sum_{\beta \in {\cal I}} P[X_\beta=1]
 E[f(Y)| X_\beta = 1].
\enas
As exchangeability implies that $E[f(Y)| X_\beta = 1]$ does not depend on $\beta$, we have
\beas
E[Y f(Y)] =
 \left( \sum_{\beta \in {\cal I}}  P[X_\beta=1] \right)
E[f(Y)| X_\alpha = 1]
= E[Y ]
E[f(Y^\alpha)],
\enas
demonstrating the first result.  The second follows easily using that $Y^I$ is a mixture of random variables all of which have distribution $Y^s$. \bbox

With $n \ge d$ we prove Theorem \ref{thm:occ} by constructing a size bias coupling of $Y_n^s$ to $Y_n$ for the urn model and verifying the hypotheses of Theorem \ref{thm:main}. To apply Lemma \ref{sblem} we construct, for each $i \in \{1,\ldots,m\}$, a configuration ${\bf M}_n^i$ that has the conditional distribution of ${\cal M}(n,m)$ given that urn $i$ contains $d$ balls on the same space as a configuration ${\bf M}_n$ with the unconditional distribution ${\cal M}(n,m)$.

We now describe the joint construction of ${\bf M}_n^i$ and ${\bf M}_n$ formally; in its course we will also define the vector ${\bf R}_n^i$ specifying the difference, up to sign, between ${\bf M}_n$ and ${\bf M}_n^i$. For a vector ${\bf M}$ and a given $i \ge 1$, let $\langle {\bf M}  \rangle_i$ be the vector obtained by deleting the $i^{th}$ component of ${\bf M}$.

With $i \in \{1,\ldots,m\}$ we first specify the $i^{th}$ components of ${\bf M}_n$ and ${\bf M}_n^i$ by letting $M_n(i) \sim {\cal B}(n,1/m)$ and $M_n^i(i)=d$, respectively. Next, let vectors ${\bf M}_{n,i}'$ and ${\bf R}_n^i$
satisfy $M_{n,i}'(i)=R_n^i(i)=0$, and whose remaining components are conditionally independent given $M_n(i)$, with conditional distributions given $M_n(i)$ specified by
$${\cal L}(\langle {\bf M}_{n,i}' \rangle_i|M_n(i)) = {\cal M}(n-M_n(i) \vee d, m-1)
$$
and
\begin{equation}
\label{cond.dist.R}
{\cal L}(\langle {\bf R}_n^i \rangle_i|M_n(i)) = {\cal M}(\abs{d-M_n(i)} , m-1),
\end{equation}
and set
\begin{equation*}
\langle {\bf M}_n \rangle_i=\langle {\bf M}_{n,i}' \rangle_i + {\bf 1}(M_n(i)<d)\,\langle {\bf R}_n^i\rangle_i \qmq{and} \langle {\bf M}_n^i \rangle_i = \langle {\bf M}_{n,i}'\rangle_i + {\bf 1}(M_n(i)>d)\,\langle {\bf R}_n^i \rangle_i.
\end{equation*}
By the additive property of the multinomial distribution, conditional on $M_n(i)$ we have that $\langle {\bf M}_n \rangle_i \sim {\cal M}(n-M_n(i),m-1)$ in all cases, so that ${\bf M}_n \sim {\cal M}(n,m)$, as required. Likewise in all cases $\langle {\bf M}_n^i \rangle_i \sim {\cal M}(n-d,m-1)$, so
\bea \label{XniXn}
{\bf M}_n \sim {\cal M}(n,m) \qmq{and} {\cal L}({\bf M}_n^i)={\cal L}({\bf M}_n|M_n(i)=d).
\ena
Further, we note that that the difference between the two configurations excluding urn $i$ satisfies
\bea \label{three.cases.szb}
\langle {\bf M}_n^i \rangle_i - \langle {\bf M}_n \rangle_i  = \mbox{sign}(M_n(i)>d)\langle {\bf R}_n^i \rangle_i, \qmq{where} \sum_{j\ge 1} R_n^i(j)=|d-M_n(i)|.
\ena

Applying the indicator function ${\bf 1}(\cdot =d)$ coordinate-wise to (\ref{XniXn}) and recalling (\ref{Yn:balls}) we obtain
\beas
{\cal L}(X_{n,1}^i,\ldots,X_{n,m}^i)={\cal L}(X_{n,1},\ldots,X_{n,m}|M_n(i)=d),
\enas
and Lemma \ref{sblem} now yields that $Y_n^i$, counting the number of urns containing $d$ balls in the configuration ${\bf M}_n^i$, given explicitly by
\beas
Y_n^i=\sum_{j\ge 1} X_{n,j}^i, \qmq{with}  X_{n,j}^i={\bf 1}(M_n^i(j)=d) \quad \mbox{for $j \ge 1$,}
\enas
has the $Y_n$-size biased distribution. Again by Lemma \ref{sblem}, if $I_n$ is uniformly distributed over $\{1,\ldots,m\}$, independent of all other variables, then $Y_n^s = Y_n^{I_n}$ also has the $Y_n$-size bias distribution.

\section{Auxiliary Results and Proof of Theorem 1.1}
\label{proof.of.theorem}

To prove Theorem~\ref{thm:occ} we utilize a general result of \citet{Go12}, given as Theorem \ref{thm:main} below, whose framework has already been described in Section~\ref{sec:intro}. In particular, the random variables of interest $Y_n, n \ge n_0$ have distributions ${\cal L}_\theta(Y_n)$ that depend on a parameter $\theta$ in a topological space $\Theta_n$, also endowed with a $\sigma$-algebra of subsets. In our application we give $\Theta_n=\mathbb{N}_2$ the discrete topology, and the $\sigma$-algebra the collection of all its subsets.

In Theorem \ref{thm:main}, $r_{n,\theta}$
is a function that determines the quality of the bound to the normal, the sequence $s_{n,\theta}$ is used to control a random variable
$L_n$  determining the size of the smaller subproblem $V_n$ related to $Y_n$. In general, the mean $\mu_{n,\theta}$ and variance $\sigma_{n,\theta}^2$ of $Y_n$ under ${\cal L}_\theta$, and $r_{n,\theta}$, are required to be measurable in $\theta$, a condition satisfied for all natural examples, and in particular, for the one considered here.
\begin{theorem}
\label{thm:main} For some $n_0 \in \mathbb{N}_0$ and all $n \ge n_0$  let
$Y_n$ be a nonnegative random variable with mean $\mu_{n,\theta}=E_\theta Y_n$ and positive variance $\sigma_{n,\theta}^2=\mbox{Var}_\theta(Y_n)$
for all $\theta \in \Theta_n$, 
and set
\bea \label{def:Wntheta}
W_{n,\theta}=\frac{Y_n-\mu_{n,\theta}}{\sigma_{n,\theta}},
\ena
the standardized value of $Y_n$. Let $r_{n,\theta}$ be positive for all $n \ge n_0$ and all $\theta \in \Theta_n$, and for all $r \ge 0$ let
\bea \label{def:Thetanr}
\Theta_{n,r} = \{\theta \in \Theta_n: r_{n,\theta} \ge r\}.
\ena
Assume there exists $r_1>0$ and $n_1 \ge n_0$ such that
\bea \label{bd.away.zero.rnt}
\max_{n_0 \le n < n_1} \sup_{\theta \in \Theta_{n,r_1}} r_{n,\theta} < \infty.
\ena
Further, suppose that for all $n \ge n_1$ and $\theta \in \Theta_{n,r_1}$, there exist random variables $Y_n^s, K_n, L_n,$ $\psi_{n,\theta}, V_n$ and $B_n$ on the same space as $Y_n$, and a $\sigma$-algebra ${\cal F}_n$, generated by a collection of random elements $J_n$, such that the following conditions hold.

\begin{enumerate}
\item \label{Psinisrootn} The random variable $Y_n^s$ has the $Y_n$-size bias distribution, and
\bea \label{def:psin}
\Psi_{n,\theta}=\sqrt{\mbox{Var}_\theta \left(E_\theta(Y_n^s-Y_n|Y_n)\right)} \qmq{satisfies} \sup_{n \ge n_1,\theta \in \Theta_{n,r_1}}\frac{r_{n,\theta}\mu_{n,\theta} \Psi_{n,\theta}}{\sigma_{n,\theta}^2} < \infty.
\ena


\item \label{assumption:Kn} The random variable $K_n$ is ${\cal F}_n$-measurable, $|Y_n^s-Y_n| \le K_n$, and
\bea \label{assumption:Kn.inequality}
\sup_{n \ge n_1,\theta \in \Theta_{n,r_1}}\frac{r_{n,\theta}\mu_{n,\theta} E_\theta[\left( 1+ |W_{n,\theta}| \right)K_n^2]}{\sigma_{n,\theta}^3}<\infty,
\ena
with $W_{n,\theta}$ as given in (\ref{def:Wntheta}).

\item \label{assumption:condlawYnvgivenK}  The random variable $L_n$ takes values in $\{0,1,\ldots,n\}$, there exists a positive integer valued sequence $\{s_{n,\theta}\}_{n \ge n_1}$ satisfying $n-s_{n,\theta} \ge n_0$, the variables $L_n$ and $\psi_{n,\theta}$ are ${\cal F}_n$-measurable, for some 
$F_{n,\theta} \in {\cal F}_n$ satisfying $F_{n,\theta} \subset \{L_n \le s_{n,\theta}\}$,
\bea \label{Vn.conditioned.Jn}
\psi_{n,\theta} \in \Theta_{n-L_n} \qmq {and} {\cal L}_\theta(V_n|J_n)={\cal L}_{\psi_{n,\theta}}(Y_{n-L_n}) \qmq{on $F_{n,\theta}$}
\ena
and
\bea \label{mu.sig.K.L.s}
\sup_{n \ge n_1,\theta \in \Theta_{n,r_1}}\frac{r_{n,\theta}^2\mu_{n,\theta}}{\sigma_{n,\theta}^3} E_\theta \left[  K_n^2 (1-{\bf 1}(F_{n,\theta}))\right] < \infty.
\ena

\item \label{unif.var.bnd}
There exist $\{\Cl{sig-ratio}, \Cl{r-ratio}\} \subset (0,\infty)$ such that
\beas 
\sigma^2_{n,\theta} \le \Cr{sig-ratio} \sigma_{n-L_n,\psi_{n,\theta}}^2
\qmq{and}
r_{n,\theta} \le \Cr{r-ratio} r_{n-L_n,\psi_{n,\theta}}
\qmq{on $F_{n,\theta}$.}
\enas

\item \label{assumption:YnmYnVnprime} The random variable $B_n$ is ${\cal F}_n$-measurable, $|Y_n-V_n| \le B_n$ and
\bea \label{inequality:YnmYnVnprime}
\sup_{n \ge n_1,\theta \in \Theta_{n,r_1}}\frac{r_{n,\theta}^2\mu_{n,\theta}
E_\theta[K_n^2 B_n]
}{\sigma_{n,\theta}^4}<\infty.
\ena

\item \label{Theta.n.r-compact} Either
\begin{enumerate}

\item  there exists $l_{n,0} \in \mathbb{N}_0$ such that $P_\theta(L_n=l_{n,0})=1$ for all $\theta \in \Theta_{n,r_1}$

or

\item the set $\Theta_{n,r_1}$ is a compact subset of $\Theta_n$, and the functions of $\theta$ 
\beas 
t_{n,\theta,l}= E_\theta \left( \frac{K_n^2}{E_\theta K_n^2}{\bf 1}(L_n=l)\right) \qm{for $l \in \{0,1,\ldots,n\}$}
\enas
are continuous on $\Theta_{n,r_1}$ for $l \in \{0,1,\ldots,s_n\}$ where $s_n=\sup_{\theta \in \Theta_{n,r_1}}s_{n,\theta}$.

\end{enumerate}

\end{enumerate}

Then there exists a constant $C$ such that for all $n \ge n_0$ and $\theta \in \Theta_n$
\beas 
\sup_{z \in \mathbb{R}}|P_\theta(W_{n,\theta}  \le z) - P(Z \le z)| \le C/r_{n,\theta}.
\enas
\end{theorem}

When higher moments exist a number of the conditions of the theorem
may be verified using standard inequalities. In particular, by the Cauchy-Schwarz inequality a sufficient condition for (\ref{assumption:Kn.inequality}) is
\bea
\label{kn-sum2-conv}
\sup_{n \ge n_1,\theta \in \Theta_{n,r_1}}\frac{r_{n,\theta}\mu_{n,\theta} k_{n,\theta,4}^{1/2}}{\sigma_{n,\theta}^3}<\infty \quad \mbox{where} \quad k_{n,\theta,m}=E_\theta K_n^m,
\ena
and, when $F_{n,\theta}=\{L_{n,\theta} \le s_{n.\theta}\}$ then a sufficient condition for (\ref{mu.sig.K.L.s}) is
\bea \label{mu.sig.K.L.s.sufficient}
\sup_{n \ge n_1,\theta \in \Theta_{n,r_1}}\frac{r_{n,\theta}^2\mu_{n,\theta}  k_{n,\theta,4}^{\frac{1}{2}}  l_{n,\theta,2}^{\frac{1}{2}}}{\sigma_{n,\theta}^3 s_{n,\theta}} < \infty \qmq{where} l_{n,\theta,m}=E_\theta L_n^m,
\ena
since, additionally using the Markov inequality yields
\beas
E_\theta \left[  K_n^2 {\bf 1}(L_n > s_{n,\theta})\right]
\le k_{n,\theta,4}^{1/2}P_\theta(L_n>s_{n,\theta})^{\frac{1}{2}}
= k_{n,\theta,4}^{1/2}P_\theta(L_n^2>s_{n,\theta}^2)^{\frac{1}{2}}
\le \frac{k_{n,\theta,4}^{\frac{1}{2}}l_{n,\theta,2}^{\frac{1}{2}}}{s_{n,\theta}}.
\enas
Similarly, a sufficient condition for (\ref{inequality:YnmYnVnprime}) is
\beas 
\sup_{n \ge n_1,\theta \in \Theta_{n,r_1}}\frac{r_{n,\theta}^2\mu_{n,\theta} k_{n,\theta,4}^{\frac{1}{2}}b_{n,\theta,2}^{\frac{1}{2}}}{\sigma_{n,\theta}^4}<\infty \quad \mbox{where} \quad b_{n,\theta,m}=E_\theta B_n^m.
\enas

 Regarding (\ref{Vn.conditioned.Jn}) we remark that by ${\cal L}_\theta(Y_{n-L_n})$ we mean the mixture distribution\\ $\sum_{m=0}^n {\cal L}_\theta(Y_m)P(L_n=n-m)$, which can be defined without requiring that $Y_0,\ldots,Y_n$ 
 and $L_n$ all be defined on the same space.

Recalling $\mathbb{N}_k = \mathbb{Z} \cap [k,\infty)$, applying Theorem \ref{thm:main} to the occupancy problem we let $n_0 = d$, $r_{n,m}$ be given by (\ref{rnm:def}), and $\Theta_n=\mathbb{N}_2$ for all $n \ge n_0$, making note of the identification between positive integers $m$ and elements given by (\ref{theta_m}) that lie in the set $\Theta$ of (\ref{def:Theta}).

Before starting the proof of Theorem \ref{thm:occ} we collect some crucial facts needed later regarding the behavior of the mean and variance of $Y_n$. Letting
\begin{equation}\label{tau}
\tau_d(x)=\frac{e^{-x}x^d}{d!} \qmq{and} \varphi_d(x)=1-\tau_d(x)-\tau_d(x)\frac{(x-d)^2}{x},
\end{equation} \citet[][p.~37-38]{Kolchin78} show that, for all $n,m\ge 1$ and $d\ge 0$,
\begin{equation}\label{muupbd}
\mu_{n,m}\le m\tau_d(n/m)e^{d/m}
\end{equation} and that for $n,m\To\infty$ such that $n/m=o(m)$,
\beas
\mu_{n,m}&=m\tau_d(n/m)+\tau_d(n/m)\left(d-\frac{n/m}{2}-\frac{{d\choose 2}}{n/m}\right)+O(1/m) 
\enas
and
\bea
\sigma_{n,m}^2&=m\tau_d(n/m)\varphi_d(n/m)(1+o(1)).\label{Kolsig}
\ena
The following lemma gives further properties of $\mu_{n,m}$ and $\sigma_{n,m}^2$, and is proved in the Appendix.

\begin{lemma} \label{r1.large}
\begin{enumerate}
\item\label{Theta.finite} For any fixed $n \ge d \ge 2$,
\bea \label{sig-to-0-m}
\lim_{m \rightarrow \infty} \sigma_{n,m}^2 = 0
\ena
and the set $\Theta_{n,r_1}$, given in (\ref{def:Thetanr}) with $r_{n,m}$ as in (\ref{rnm:def}), is finite for all $r_1>0$.

\item\label{lem:alln,m} Let $d\ge 1$. There are constants $\Cl{sigma.over.mu}, \Cl{mu<n}, \Cl{sig<n}$, depending only on $d$, such that, for all $n\ge d$ and $m\ge 2$,
\begin{enumerate}
\item \label{sig2.le.mu}  $\sigma_{n,m}^2 \le \Cr{sigma.over.mu} \mu_{n,m}$
\item\label{mu,sig<n} $\mu_{n,m}\le \Cr{mu<n}n$  and $\sigma_{n,m}^2\le \Cr{sig<n} n$.
\end{enumerate}

\item\label{infvarphi>0} Let $d\ge 2$. With $\varphi_d(x)$ given by \eqref{tau}, $\inf_{x>0}\varphi_d(x)>0$.

\item\label{lem:sig>r1} Let $d\ge 2$. Given any $n^*, m^*$ and $\eps>0$ there are constants $r_1$ and $\Cl{sig/mu>C}$ such that all $n,m$ satisfying $\sigma_{n,m}^2\ge r_1$ also satisfy
\begin{enumerate}
\item\label{nmToinf} $n>n^*$ and $m>m^*$
\item\label{n/mlogm} $n/m\le (1+\eps)\log m$
\item\label{sig/mu>C} $\mu_{n,m}\le \Cr{sig/mu>C}\sigma_{n,m}^2$.
\end{enumerate}

\end{enumerate}
\end{lemma}

In keeping with the notation of Theorem \ref{thm:main} and the identification between elements in $\mathbb{N}_1$ and $\Theta$ in (\ref{def:Theta}) and as described at the end of Section \ref{sec:intro}, in the following we will use $E_m$, $\var_m$, and $P_m$ to respectively denote expectation, variance, and probability with respect to a multinomial distribution with probability parameter \eqref{theta_m}.

\subsection*{Proof of Theorem 1.1}
We prove Theorem~\ref{thm:occ} by verifying the conditions of Theorem \ref{thm:main}. When $n=d$ and $m \ge 2$, the probability that all $d$ balls fall in urn 1 is positive, as is the probability that $d-1$ balls fall in urn 1. Hence $P(Y_n>0)$ and $P(Y_n=0)$ are both positive, so $Y_n$ is not constant almost surely, and its variance $\sigma_{d,m}^2$ is strictly positive. The same conclusion holds for $n \ge d+1$ and $m \ge 2$ by considering the event that $d$ balls fall in urn 1 and $n-d$ in urn 2, and the event that all balls fall in urn 1. Hence $r_{n,m}$ given in (\ref{rnm:def}) is also positive for all $n \ge n_0$ and $m \in \Theta_n$. In lieu of naming $n_1$ and $r_1$ explicitly, we show that the conditions of Theorem \ref{thm:main} are satisfied by choosing $n_1$ and $r_1$ sufficiently large.  By Part~\ref{Theta.finite} of Lemma~\ref{r1.large}, the set $\Theta_{n,r_1}$ is finite for all $n \ge n_0$ and $r_1>0$, hence (\ref{bd.away.zero.rnt}) is satisfied for any such pair.

To help with the verification of the six conditions of Theorem \ref{thm:main} we first note that Parts~\ref{lem:alln,m} and \ref{lem:sig>r1} of Lemma~\ref{r1.large} allow us to choose $r_1>0$ such that there exist positive constants $\Cr{sigma.over.mu}, \Cr{sig<n}, \Cr{sig/mu>C}$ such that $\sigma_{n,m}^2\ge r_1$ implies
\begin{gather}
\sigma_{n,m}^2 \le \Cr{sigma.over.mu} \mu_{n,m}\label{pf:sig<Cmu}\\
\sigma_{n,m}^2\le  \Cr{sig<n} n\label{r1:sig^2<Cn}\\
\mu_{n,m}\le \Cr{sig/mu>C}\sigma_{n,m}^2\label{pf:mu<Csig}\\
\frac{n}{m}\le 2\log m \le\sqrt{m/3}\label{n/m<log<sqrtm}\\
m\ge 3.\label{m>=3}
\end{gather}
Below we will repeatedly use these bounds along with the fact that
\begin{equation}
\label{r.min.subset}
\Theta_{n,r_1} \subset \left\{ m: \sigma_{n,m}^2 \ge r_1 \right\} \qmq{for all $n \ge d$ and $m \ge 2$,}
\end{equation}
which follows from directly from \eqref{def:Thetanr} and the fact that $r_{n,m} \le \sigma_{n,m}$.

\subsubsection*{Verification of Condition \ref{Psinisrootn}}
We verify inequality \eqref{def:psin} in Condition \ref{Psinisrootn} of Theorem \ref{thm:main} by showing that, for $Y_n^s$ constructed as in Section~\ref{size.bias.coupling}, there is a constant~$\Cl{Psinm}$ and integer $n_1 \in \mathbb{N}_1$ such that for all $n \ge n_1$ and $m \in \Theta_{n,r_1}$,
the quantity $\Psi_{n,m}$ satisfies
\bea \label{Psi.bound}
\Psi_{n,m}  \le \Cr{Psinm}\frac{1+(n/m)^3}{\sqrt{n}}.
\ena
Inequality (\ref{Psi.bound}) implies \eqref{def:psin} as
$$\frac{r_{n,m}\mu_{n,m} \Psi_{n,m}}{\sigma_{n,m}^2} \le \Cr{sig/mu>C}r_{n,m} \Psi_{n,m}
\le
 \Cr{sig/mu>C} \Cr{Psinm}\frac{\sigma_{n,m}}{\sqrt{n}} \le \Cr{sig/mu>C} \Cr{Psinm} \sqrt{\Cr{sig<n}},$$ where we have used \eqref{pf:mu<Csig} and \eqref{r1:sig^2<Cn}. Hence we turn our attention to showing (\ref{Psi.bound}).

By conditional Jensen's inequality, as $Y_n$ is a function of ${\bf M}_n$,
\beas
\var_m(E_m(Y_n^s-Y_n|Y_n)) \le \var_m(E_m(Y_n^s-Y_n|{\bf M}_n)).
\enas
Recalling that $I_n$ is chosen uniformly from $\{1,\ldots,m\}$, independently of the configuration ${\bf M}_n$, and that $X_j^{I_n}, j \ge 1$ is the indicator that urn $j$ contains exactly $d$ balls in the size biased configuration, we have that
$$
Y_n^s-Y_n=\sum_{j\ge 1} (X_{n,j}^{I_n}-X_{n,j})=(X_{n,I_n}^{I_n}-X_{n,I_n})+\sum_{j \not = I_n} (X_{n,j}^{I_n}-X_{n,j}).
$$
Averaging over $I_n$, we obtain
\bea
E_m[Y_n^s-Y_n|{\bf M}_n] &=&
\frac{1}{m}\sum_{i=1}^m {\bf 1}(M_n(i) \not = d) +
\frac{1}{m} \sum_{1 \le i,j \le m, j \not =i} P_m(M_n^i(j)=d|{\bf M}_n){\bf 1}(M_n(j) \not  =d)\nn \\
&&  - \frac{1}{m} \sum_{1 \le i,j \le m, j \not =i}P_m(M_n^i(j)\not =d|{\bf M}_n){\bf 1}(M_n(j)=d). \label{L-MV:27}
\ena
To understand the first sum, note that since urn $I_n$ always contains $d$ balls in the size biased configuration, $X_{n,I_n}^{I_n}-X_{n,I_n} = 1 - {\bf 1}(M_n(I_n)=d))={\bf 1}(M_n(I_n)\not =d))$, so averaging over $I_n$, which takes the values $1,\ldots,m$ each with probability $1/m$, yields the first term. The next two terms arise from the fact that $X_{n,j}^{I_n}-X_{n,j} \in \{-1,0,1\}$; in particular, the second term accounts for the cases when this difference is $1$, and the third term for when it is $-1$. For the second sum, when $I_n=i$ we have $X_{n,j}^i-X_{n,j}=1$ for $j \not =i$ if and only if $X_{n,j}^i=1$ and $X_{n,j}=0$, that is, if and only if $M_n^i(j)=d$ and $M_n(j)\not =d$. Likewise, for the third sum,
$X_{n,j}^i-X_{n,j}=-1$ for $j \not =i$ if and only if $X_{n,j}^i=0$ and $X_{n,j}=1$, and so if and only if $M_n^i(j)\not =d$ and $M_n(j)=d$.

To obtain a bound on the variance of $E_m[Y_n^s-Y_n|{\bf M}_n]$ we apply the inequality
\bea \label{Vark}
\var\left(\sum_{i=1}^k A_i\right) \le k\sum_{i=1}^k \var(A_i)
\ena
in order to handle the terms arising from (\ref{L-MV:27}) separately. We will use \eqref{Vark} and $(\sum_{i=1}^k c_i )^2\le k\sum_{i=1}^k c_i^2$ for any $c_1,\ldots,c_k$,  repeatedly below without further mention. The factor of $1/m$ outside each sum in \eqref{L-MV:27} contributes a factor of $1/m^2$ to the variance, which is withheld until further notice below.

To bound the variance of the first sum, we note that
\begin{multline}
\var_m\left( \sum_{i=1}^m {\bf 1}(M_n(i) \not = d) \right)=\var_m\left( m-\sum_{i=1}^m {\bf 1}(M_n(i) = d) \right)
=\var_m(m-Y_n)=\var_m(Y_n)\\
=\sigma_{n,m}^2\le \Cr{sig<n}n\label{Ass1-1stsum}
\end{multline} by \eqref{r1:sig^2<Cn}.

For considering the calculation of the variance for the next sum, as $M_n(j)=M_n^i(j)$ when $M_n(i)=d$ by (\ref{three.cases.szb}), we have
$$
P_m(M_n^i(j)=d|{\bf M}_n){\bf 1}(M_n(j) \not  =d)=P_m(M_n^i(j)=d|{\bf M}_n){\bf 1}(M_n(i) \not = d, M_n(j) \not  =d),
$$
and therefore may write
\beas 
\sum_{1 \le i,j \le m, i \not = j} P_m(M_n^i(j)=d|{\bf M}_n){\bf 1}(M_n(j) \not  =d) =\sum_{1 \le i,j \le m,i \not = j} a_{n,m}(i,j) + \sum_{1 \le i,j \le m,i \not = j} b_{n,m}(i,j)
\enas
where for $i \not = j$ we set
\bea \label{ab:defs}
a_{n,m}(i,j)&=& P_m(M_n^i(j)=d|{\bf M}_n){\bf 1}(M_n(i)>d,M_n(j) \not  =d),\\
b_{n,m}(i,j)&=& P_m(M_n^i(j)=d|{\bf M}_n){\bf 1}(M_n(i)<d, M_n(j)\not  =d). \nonumber
\ena
For considering the third sum in \eqref{L-MV:27}, we also define
\bea \label{cnmij:def}
c_{n,m}(i,j)=P_m(M_n^i(j)\not =d|{\bf M}_n){\bf 1}(M_n(j)=d).
\ena
In Lemma~\ref{lem:E-Sabc}, following the proof of this theorem,  it is shown that there are constants $\Cl{varsump}, \Cl{varsumb}$ and $\Cl{Ass1-3rdsum}$ such that
\begin{align}
\var_m\left(\sum_{1 \le i,j \le m,i\ne j}a_{n,m}(i,j)\right) &\le \Cr{varsump} n \left[1+\left(\frac{n}{m}\right)^4\right]\label{varsump}\\
\var_m\left(\sum_{1 \le i,j \le m,i\ne j}b_{n,m}(i,j)\right)&\le \Cr{varsumb} \frac{m^2}{n}\label{varsumb}\\
\var_m\left(\sum_{1 \le i,j \le m,i\ne j}c_{n,m}(i,j)\right)&\le \Cr{Ass1-3rdsum} n\left[1+\left(\frac{n}{m}\right)^2\right]\label{Ass1-3rdsum}
\end{align} for all $n \ge n_1$ and $m \in \Theta_{n,r_1}$. Combining \eqref{Ass1-1stsum} and \eqref{varsump}-\eqref{Ass1-3rdsum}, and accounting for the $1/m$ factors in \eqref{L-MV:27}, we have
\begin{align*}
\Psi_{n,m}^2&\le \frac{4}{m^2}\C \left\{n+ n\left[1 +\left(\frac{n}{m}\right)^4 \right]+\frac{m^2}{n}+n\left[1 +\left(\frac{n}{m}\right)^2 \right]  \right\}\\
 &\le \Cl{Psi^2} \frac{n}{m^2} \left\{1+\left(\frac{m}{n}\right)^2 +\left(\frac{n}{m}\right)^2 +\left(\frac{n}{m}\right)^4\right\}\\
 &= \Cr{Psi^2} \frac{1}{n} \left\{\left(\frac{n}{m}\right)^2+1 +\left(\frac{n}{m}\right)^4 +\left(\frac{n}{m}\right)^6\right\}\\
 &\le 3\Cr{Psi^2} \frac{1}{n} \left\{1 +\left(\frac{n}{m}\right)^6\right\},
 \end{align*} where in the last step and below we use the elementary bound
 \begin{equation}\label{xpoly-bd}
1+x^{\ell_1}+x^{\ell_2}\ldots+x^{\ell_j}\le (j+{\bf 1}\{\ell_j<\ell\})(1+x^\ell)\qmq{for all}x>0,\; 1\le\ell_1\le\ldots \le\ell_j\le\ell.
\end{equation}
Then taking $\Cr{Psinm}=\sqrt{3\Cr{Psi^2}}$ yields
$$\Psi_{n,m}\le \Cr{Psinm}\frac{\sqrt{1+(n/m)^6}}{\sqrt{n}}\le \Cr{Psinm}\frac{1+(n/m)^3}{\sqrt{n}}.$$

\subsubsection*{Verification of Condition \ref{assumption:Kn}: $K_n$ and its moments}

Let $J_n=(I_n, M_n(I_n))$, the ordered pair consisting of the identity $I_n$ of the selected urn and the number $M_n(I_n)$ of balls it contains, and recall that ${\cal F}_n$ is the $\sigma$-algebra generated by $J_n$. For $D \sim {\cal B}(n,p)$ and $q \in \mathbb{N}_1$ we have
\bea \label{bin:mom}
E D^q = \sum_{j=1}^q S_{j,q} (n)_j p^j \le \sum_{j=1}^q S_{j,q} n^j p^j \le \Cl{Riordin1}{_{,q}}(np+(np)^q)
\le \Cl{Riordin}{_{,q}} (1+(np)^q),
\ena
where in the first equality, due to \citet{Riordan37}, $S_{j,q}$ are the Stirling numbers of the second kind and $(n)_j$ is the falling factorial, and in the second inequality $\Cr{Riordin1}{_{,q}}  =   q\max_{1 \le j \le q}S_{j,q}$.

Clearly
\beas 
K_n=1+|d-M_n(I_n)|
\enas
is ${\cal F}_n$-measurable, being a function of $M_n(I_n)$. Recalling (\ref{three.cases.szb}) from the construction
in Section~\ref{size.bias.coupling}, accounting for urn $I_n$ we see that the occupancy of at most $K_n$ urns are different in the configurations ${\bf M}_n^i$ and ${\bf M}_n$ for any $i$. In particular, $|Y_n^s - Y_n| \le K_n$.

By the triangle inequality $K_n \le (1+d) + M_n(I_n)$, and taking $q^{th}$ power, by a standard inequality and (\ref{bin:mom}) we obtain
\begin{multline}
E_mK_n^q \le 2^{q-1} \left( (1+d)^q + E_mM_n(I_n)^q \right) \le
2^{q-1}  \left( (1+d)^q + \Cr{Riordin}{_{,q}}\left(1 + (n/m)^q \right) \right) \\
\le \Cl{Kn1.moment}{_{,q}}\left( d^q + 1 + (n/m)^q \right)
\le \Cl{Kn.moment}{_{,q}}\left( 1 +(n/m)^q \right). \label{Kn.moment.bound.q}
\end{multline}

We now show that (\ref{kn-sum2-conv}), sufficient for (\ref{assumption:Kn.inequality}), is satisfied. Applying the definition (\ref{rnm:def}) of $r_{n,m}$, \eqref{pf:mu<Csig}, \eqref{r.min.subset}, and the moment bound (\ref{Kn.moment.bound.q}),  there is some $n_1$ such that for all $n \ge n_1$ and any $m \in \Theta_{n,r_1}$,
\beas
\frac{r_{n,m}\mu_{n,m} k_{n,m,4}^{1/2}}{\sigma_{n,m}^3}
\le \Cr{sig/mu>C}\frac{\left[\Cr{Kn.moment}{_{,2}}(1+(n/m)^4)\right]^{1/2}}{1+(n/m)^3}\le \Cr{sig/mu>C}\sqrt{\Cr{Kn.moment}{_{,2}}}\left(\frac{ 1+(n/m)^2}{1+(n/m)^3}\right)\le 2\Cr{sig/mu>C}\sqrt{\Cr{Kn.moment}{_{,2}}},
\enas using \eqref{xpoly-bd} in this last step.

\subsubsection*{Verification of Condition \ref{assumption:condlawYnvgivenK}: $L_n$ and its moments}
Set
\beas
L_n=M_n(I_n), \quad \psi_{n,m}=m-1, \quad s_{n,m}=\left\lceil n^{1/2}\, \right\rceil \qmq{and} F_{n,m}=\{L_n \le s_{n,m}\}.
\enas
Clearly $L_n$ takes values in $\{0,1,\ldots,n\}$, and $n-s_{n,m} \ge n_0$ for all $n$ sufficiently large,  and $L_n, \psi_{n,m}$ and $F_{n,m}$ are ${\cal F}_n$ measurable.  Now, by \eqref{m>=3},  the first part of (\ref{Vn.conditioned.Jn}) holds.

Let
\bea \label{Vn:def}
V_n=\sum_{i \not = I_n} X_{n,i}
\ena
with $X_{n,i}$ as in (\ref{Yn:balls}). Using that $I_n$ is independent of $M_n(j),j=1,\ldots,m$, and the properties of the multinomial ${\cal M}(n,m)$ distribution, we have
\beas 
{\cal L}(M_n(j), j \not =I_n|M_n(I_n)=l,I_n=i) = {\cal L}(M_n(j), j \not =m|M_n(m)=l)= {\cal M}(n-l,1/(m-1)),
\enas
and counting the number of urns with occupancy $d$ on both sides of this distributional identity yields
\beas
{\cal L}_m(V_n|J_n)={\cal L}_{m-1}(Y_{n-M_n(I_n)}) = {\cal L}_{\psi_{n,m}}(Y_{n-L_n}),
\enas
so the second part of (\ref{Vn.conditioned.Jn}) holds on the entire probability space, so in particular on $F_{n,\theta}$. As $L_n \sim {\cal B}(n,1/m)$ under $P_m$, from (\ref{bin:mom}) we obtain
\bea \label{Ln.moment.bound.q}
E_m L_n^q  \le \Cr{Riordin}{_{,q}} (1+(n/m)^q).
\ena
Hence, inequality (\ref{mu.sig.K.L.s.sufficient}), sufficient for (\ref{mu.sig.K.L.s}), holds as
\begin{multline*}
\frac{r_{n,m}^2\mu_{n,m}  k_{n,m,4}^{\frac{1}{2}}  l_{n,m,2}^{\frac{1}{2}}}{\sigma_{n,m}^3 s_{n,m}}
\le \C\frac{\mu_{n,m}\sqrt{1+(n/m)^4}\sqrt{1+(n/m)^2}}{\sigma_{n,m} \sqrt{n} \left[1 + (n/m)^3\right]^2}
\le \Cl{L.2}\frac{\sigma_{n,m}\sqrt{1+(n/m)^6}}{\sqrt{n}\left[1 + (n/m)^6\right]} \\
\le \Cr{L.2}\sqrt{\Cr{sig<n}},
\end{multline*}
where we have used the definition (\ref{rnm:def}) of $r_{n,m}$, the definition of $s_{n,m}$, (\ref{Kn.moment.bound.q}) and (\ref{Ln.moment.bound.q}) in the first inequality, \eqref{pf:mu<Csig} in the second inequality, and \eqref{r1:sig^2<Cn} in the final inequality.

\subsubsection*{Verification of Condition \ref{unif.var.bnd}}

We first show that there exists $n_1$ such that, for all $n \ge n_1$ and $m \in \Theta_{n,r_1}$,
\bea \label{mu.n.m.n-L}
\mu_{n,m} \le 18 \mu_{n-L_n, m-1} \qm{on $F_{n,m}$.}
\ena
As $n/(n-\lceil \sqrt{n} \rceil) \rightarrow 1$ as $n \rightarrow \infty$ and $F_{n,m}=\{L_n \le \lceil \sqrt{n} \rceil\}$, there exists $n_1$ such that $n-\lceil \sqrt{n} \rceil \ge n_0$ and
\beas
\frac{(n)_d}{(n-L_n)_d} &\le& \left( \frac{n}{n-\lceil \sqrt{n} \rceil} \right)^d \le 2 \qmq{on $F_{n,m}$,}
\enas
for all $n \ge n_1$.

Next, as $m \ge 3$ by \eqref{m>=3} we obtain $m^2-2m \ge m^2/3$, and therefore, using the first upper bound on $n/m$ in (\ref{n/m<log<sqrtm}) for the second to last inequality, we obtain
\beas
\frac{\left(1-\frac{1}{m}\right)^{n-d}}{\left(1-\frac{1}{m-1}\right)^{n-d}} = \left( 1+\frac{1}{m^2-2m}\right)^{n-d}
\le \left( 1+\frac{3}{m^2}\right)^n \le e^{3n/m^2} \le e^{6 \log m/m} \le 9.
\enas
Hence, for all $n \ge n_2$ and $m \in \Theta_{n,r_1}$, on $F_{n,m}$, recalling (\ref{occ:mv}) we have
\beas
\mu_{n,m} &=& \frac{(n)_d}{d!} \frac{1}{m^{d-1}}\left( 1-\frac{1}{m} \right)^{n-d}\\
&<& \frac{(n)_d}{d!} \frac{1}{(m-1)^{d-1}}\left( 1-\frac{1}{m} \right)^{n-d}\\
&\le& 2 \frac{(n-L_n)_d}{d!} \frac{1}{(m-1)^{d-1}}\left( 1-\frac{1}{m} \right)^{n-d}\\
&\le& 18 \frac{(n-L_n)_d}{d!} \frac{1}{(m-1)^{d-1}}\left( 1-\frac{1}{m-1} \right)^{n-d}\\
&\le& 18 \frac{(n-L_n)_d}{d!} \frac{1}{(m-1)^{d-1}}\left( 1-\frac{1}{m-1} \right)^{n-L_n-d}\\
&=& 18 \mu_{n-L_n,m-1}.
\enas
By \eqref{r.min.subset}, \eqref{pf:mu<Csig} and \eqref{pf:sig<Cmu} hold whenever $m \in \Theta_{n,r_1}$.
Now the first part of Condition \ref{unif.var.bnd} follows from (\ref{mu.n.m.n-L}), \eqref{pf:mu<Csig}, and \eqref{pf:sig<Cmu} since, for all $n \ge n_1$ and $m \in \Theta_{n,r_1}$,
\bea \label{unif.var.bnd.sigma}
\sigma_{n,m}^2 \le \Cr{sigma.over.mu} \mu_{n,m} \le 18 \Cr{sigma.over.mu} \mu_{n-L_n,m-1} \le \Cr{sig-ratio} \sigma_{n-L_n,m-1}^2  \qmq{on $F_{n,m}$, where $\Cr{sig-ratio}=18 \Cr{sig/mu>C}\Cr{sigma.over.mu}$.}
\ena
Since for $m\ge 2$,
$$1+\left(\frac{n}{m-1}\right)^3=1+\left(\frac{n}{m}\right)^3\left(\frac{m}{m-1}\right)^3\le 1+8\left(\frac{n}{m}\right)^3 \le 8\left[1+\left(\frac{n}{m}\right)^3\right],$$
and now the second part of Condition \ref{unif.var.bnd} follows with the help of (\ref{unif.var.bnd.sigma}) since
\beas
r_{n,m} = \frac{\sigma_{n,m}}{1+(\frac{n}{m})^3}
\le \frac{\sqrt{\Cr{sig-ratio}} \sigma_{n-L_n,m-1}}{1+(\frac{n}{m})^3}
\le \frac{8\sqrt{\Cr{sig-ratio}} \sigma_{n-L_n,m-1}}{1+(\frac{n}{m-1})^3}
\le \frac{8\sqrt{\Cr{sig-ratio}} \sigma_{n-L_n,m-1}}{1+(\frac{n-L_n}{m-1})^3}
= 8\sqrt{\Cr{sig-ratio}}r_{n-L_n,m-1}.
\enas

\subsubsection*{Verification of Condition \ref{assumption:YnmYnVnprime}: $B_n$ and its moments}
With $V_n$ given by (\ref{Vn:def}), we have $|Y_n-V_n| = X_{n,I_n} \le 1$,
so we take $B_n=1$, which is trivially ${\cal F}_n$-measurable. Now using \eqref{pf:mu<Csig}, \eqref{xpoly-bd}, and \eqref{Kn.moment.bound.q} we obtain
\beas
\frac{r_{n,\theta}^2\mu_{n,\theta}
E_m[K_n^2 B_n]
}{\sigma_{n,\theta}^4} = \frac{r_{n,\theta}^2\mu_{n,\theta}k_{n,m,2}}{\sigma_{n,\theta}^4} \le \Cr{Kn.moment}{_{,2}}
\frac{\mu_{n,m}\left(1 + (n/m)^2 \right)}{\sigma_{n,m}^2 \left(1 + (n/m)^6 \right)} \le 2\Cr{Kn.moment}{_{,2}} \Cr{sig/mu>C}.
\enas

\subsubsection*{Verification of Condition \ref{Theta.n.r-compact}}
Endowing the set $\mathbb{N}_2$ of integers with the discrete topology, a subset of $\Theta_{n,r_1} \subset \mathbb{N}_2$ is compact if and only if it is finite. As any function on a set with the discrete topology is continuous, Condition \ref{Theta.n.r-compact}b is a consequence of Lemma~\ref{r1.large}, Part~\ref{Theta.finite}. \bbox

\bigskip

Next we state and prove a lemma used in the verification of Condition \ref{Psinisrootn}.

\begin{lemma}\label{lem:E-Sabc}
Let $d \in \{2,3,\ldots\}$. There exists $n_1 \ge n_0$ and constants $\Cr{varsump}, \Cr{varsumb}, \Cr{Ass1-3rdsum}$ depending only on $d$ such that \eqref{varsump}-\eqref{Ass1-3rdsum} hold for all  $n\ge n_1$ and $m\in\Theta_{n,r_1}$.
\end{lemma}

\proof Consider first \eqref{varsump}. By (\ref{three.cases.szb}), for all $1 \le i,j \le m, i \not =j$, on $M_n(i)>d$ we have
\beas
M_n^i(j)=M_n(j)+R_n^i(j),
\enas
so that $R_n^i(j)$ is the number  of the `excess'  $M_n(i)-d$ balls distributed to urn $j$, which requires $d-M_n(j)$ of them to achieve $M_n^i(j)=d$.  Thus $a_{n,m}(i,j)=0$ unless $M_n(i)-d \ge d-M_n(j)$, that is, unless $M_n(i)+M_n(j) \ge 2d$. Hence, from (\ref{ab:defs}),
\begin{multline*}
a_{n,m}(i,j) = P_m(M_n^i(j)=d|{\bf M}_n){\bf 1}(M_n(i)>d,M_n(j) \not = d)\\
             = P_m(M_n(j)+R_n^i(j)=d|{\bf M}_n){\bf 1}(M_n(i)+M_n(j) \ge 2d,M_n(i)>d,M_n(j) \not = d)\\
             = P_m(R_n^i(j)=d-M_n(j)|{\bf M}_n){\bf 1}(M_n(i)+M_n(j) \ge 2d,M_n(i)>d,M_n(j)< d),
\end{multline*}
where we have used that $R_n^i(j) \ge 0$ makes $M_n(j) > d$ impossible in the second equality.
As $M_n(i)+M_n(j) \ge 2d$ and $M_n(j)<d$ imply that $M_n(i)>d$, letting $p=1/(m-1)$ we have, that
\begin{multline}
a_{n,m}(i,j)= P_m(R_n^i(j)=d-M_n(j)|{\bf M}_n){\bf 1}(M_n(i)+M_n(j) \ge 2d, M_n(j) < d)\\
= P_m(R_n^i(j)=d-M_n(j)|M_n(i)){\bf 1}(M_n(i)+M_n(j) \ge 2d, M_n(j) < d)\\
=\label{def.pij} {M_n(i)-d \choose d-M_n(j)} p^{d-M_n(j)}\left(1-p\right)^{M_n(i)+M_n(j)-2d}{\bf 1}(M_n(i)+M_n(j) \ge 2d, M_n(j) < d),
\end{multline}
where we have used that ${\bf M}_n$ and ${\bf R}_n^i$ are conditionally independent given $M_n(i)$, and therefore that the conditional distribution of ${\bf R}_n^i$ given ${\bf M}_n$ is the same as that given $M_n(i)$, specified in (\ref{cond.dist.R}).

Now considering $b_{n,m}(i,j)$ in (\ref{ab:defs}), using \eqref{three.cases.szb} and arguing similarly we obtain
\beas
b_{n,m}(i,j) &=& P_m(M_n^i(j)=d|{\bf M}_n){\bf 1}(M_n(i)<d, M_n(j) \not  =d)\\
&=& P_m(M_n(j)-R_n^i(j)=d|{\bf M}_n){\bf 1}(M_n(i)<d, M_n(j) \not  =d)\\
&=& P_m(R_n^i(j)=M_n(j)-d|{\bf M}_n){\bf 1}(M_n(i)<d, M_n(j) > d).
\enas
By (\ref{cond.dist.R}), $R_n^i(j) \le d-M_n(i)$ when $M_n(i)<d$, and therefore $b_{n,m}(i,j)=0$ unless $d-M_n(i)\ge M_n(j)-d$, that is, unless $M_n(i)+M_n(j)\le 2d$. Hence
\bea
b_{n,m}(i,j)  &=& P_m(R_n^i(j)=M_n(j)-d|{\bf M}_n){\bf 1}(M_n(i)+M_n(j) \le 2d, M_n(i)<d, M_n(j) > d)\nonumber \\
&=& P_m(R_n^i(j)=M_n(j)-d|{\bf M}_n){\bf 1}(M_n(i)+M_n(j) \le 2d, M_n(j) > d)\nonumber \\
&=&
\frac{{n-d\choose d}{d-M_n(i)\choose M_n(j)-d}}{{n-M_n(i)\choose M_n(j)}}{\bf 1}(M_n(i)+M_n(j)\le 2d, M_n(j) > d),
\label{bij.def1}
\ena
using that the conditional distribution of $R_n^i(j)$ given ${\bf M}_n$, as $M_n(j)=M_n^i(j)+R_n^i(j)$, is hypergeometric. As $n<2d$ implies ${n-d \choose d}=0$, which implies $b_{n,m}(i,j)=0$, we assume $n \ge 2d$ when proving (\ref{varsumb}).

Considering $c_{n,m}(i,j)$ in \eqref{cnmij:def}, let $q=1-p$ and write
\begin{multline} \label{cnm:def}
c_{n,m}(i,j) = P_m(R_n^i(j) \not =0|{\bf M}_n){\bf 1}(M_n(j)=d)=(1-P_m(R_n^i(j)=0|{\bf M}_n)){\bf 1}(M_n(j)=d)\\
=\left(1-q^{|M_n(i)-d|}\right){\bf 1}(M_n(j)=d).
\end{multline}

To prove each of \eqref{varsump}-\eqref{Ass1-3rdsum} we apply the inequality of \citet{EfSt81}. Let $S_{n-1}(x_1,\ldots,x_{n-1})$
be a symmetric function of $x_1,\ldots,x_{n-1}$, and suppose that $U_1,\ldots,U_n$ are i.i.d. random variables. For $k=1,\ldots,n$, let $S_{n,(k)}$ be the value of $S_{n-1}$ computed by omitting the $k^{th}$ variable~$U_k$, that is,
$$S_{n,(k)}=S_{n-1}(U_1,\ldots,U_{k-1},U_{k+1},\ldots,U_{n}), \qmq{and set} S_{n,(\cdot)}= \frac{1}{n}\sum_{k=1}^{n} S_{n,(k)}.$$ Then  by \citet[][Equation~1.6]{EfSt81},
\beas
\var(S_{n,(n)}) \le E \sum_{k=1}^{n} (S_{n,(k)}-S_{n,(\cdot)})^2.
\enas
As the average $S_{n,(\cdot)}$ minimizes the sum of squares, replacing it by any symmetric function $T_n$ of $U_1,\ldots,U_n$ yields
\bea \label{efron.stein.n+1.a}
\var(S_{n,(n)}) \le E \sum_{k=1}^{n} (S_{n,(k)}-T_n)^2=nE(S_{n,(n)}-T_n)^2,
\ena
this last equality since the distribution of $S_{n,(k)}-T_{n}$ does not depend on $k$.

In order to apply (\ref{efron.stein.n+1.a}), independently label the $n$ balls $1$ through $n$ such that all orderings are equally likely, and let the variables $U_k \in \{1,\ldots,m\}$ denote the location of the $k^{th}$ ball, $k=1,\ldots,n$.  Note that the three functions (\ref{def.pij}), (\ref{bij.def1}) and (\ref{cnm:def}) can be written for all $n$ as $T_n=T(n,M_n(i),M_n(j))$ for some function $T$. Hence, applying the Efron-Stein inequality with $S_{n-1}(U_1,\ldots,U_{n-1})=T(n-1,M_{n-1}(i),M_{n-1}(j))$, we obtain $S_{n,(n)}=T_{n-1}$ and \eqref{efron.stein.n+1.a} yields
\bea \label{efron.stein.n+1}
\var(T_{n-1}) \le nE(T_{n-1}-T_n)^2.
\ena
In particular, to prove \eqref{varsump} we apply (\ref{efron.stein.n+1}) with
\beas
T_n=\sum_{1 \le i,j \le m,i\ne j}a_{n,m}(i,j)\qmq{and}
T_{n-1}=\sum_{1 \le i,j \le m,i\ne j}a_{n,m,(n)}(i,j)
\enas
where $a_{n,m,(n)}(i,j)$ is the value of $a_{n,m}(i,j)$ in (\ref{def.pij})
when withholding ball $n$. As ${\cal L}(M_n^i(j)|{\bf M}_n)={\cal L}(M_n^ i(j)|M_n(i),M_n(j))$, we have that $a_{n,m,(n)}(i,j)=a_{n,m}(i,j)$  whenever $U_{n}\not\in\{i,j\}$, and hence
\begin{equation}
T_{n-1}-T_n=\sum_{1 \le j \le m, j\ne U_{n}}[a_{n,m,(n)}(U_{n},j)-a_{n,m}(U_{n},j)]+\sum_{1 \le i \le m, i\ne U_{n}}[a_{n,m,(n)}(i,U_{n})-a_{n,m}(i,U_{n})].\label{efst_mix_sums}
\end{equation}
By \eqref{def.pij} we can further restrict the summation of the first sum in (\ref{efst_mix_sums}) over indices $j$ in the union of the random index sets
\begin{align*}
J_1&=\{1 \le j \le m, j\ne U_{n}: M_n(U_n)+M_n(j)\ge 2d+1, M_n(j)<d\}\qm{and}\\
J_2&=\{1 \le j \le m, j\ne U_n: M_n(U_n)+M_n(j)= 2d, M_n(j)<d\}.
\end{align*}
For $j\in J_1$,
\begin{multline*}
a_{n,m,(n)}(U_{n},j)-a_{n,m}(U_{n},j)={M_n(U_n)-d-1\choose d-M_n(j)}p^{d-M_n(j)} (1-p)^{M_n(U_n)+M_n(j)-2d-1}\\
-{M_n(U_n)-d\choose d-M_n(j)}p^{d-M_n(j)} (1-p)^{M_n(U_n)+M_n(j)-2d}\\
={M_n(U_n)-d\choose d-M_n(j)} p^{d-M_n(j)} (1-p)^{M_n(U_n)+M_n(j)-2d-1}\left(p-\frac{d-M_n(j)}{M_n(U_n)-d}\right),
\end{multline*} and for $j \in J_1$ this last term is bounded above in absolute value by
\begin{equation}\label{absbdJ1}
\abs{p-\frac{d-M_n(j)}{M_n(U_n)-d}}\le p+\left[\max_{ x+y\ge 2d+1,\; x\le d-1}\left(\frac{d-x}{y-d}\right)\right]\le p+\frac{d-0}{(d+2)-d}\le 1+d/2=:\Cl{absbdJ1},
\end{equation} and hence
\bea \label{abs.diff.an}
|a_{n,m,(n)}(U_{n},j)-a_{n,m}(U_{n},j)| \le \Cr{absbdJ1}{M_n(U_n)-d\choose d-M_n(j)} p^{d-M_n(j)}.
\ena
To bound the right hand side we will use the fact that, for any $x, k, \ell\in \mathbb{N}_0$ satisfying $0\le x\le k-\ell$,
\begin{equation}\label{(kx)p^x-bd}
{k\choose x+\ell}p^x\le \begin{cases}
{k\choose \ell},&k\le (\ell+1)/p+\ell\\
2^k,&\mbox{for all $k$}.
\end{cases}\end{equation} The second case is trivial, and to prove the first write
$$
\frac{{k\choose x+\ell}p^x}{{k\choose \ell}}=\prod_{1\le i\le x}\frac{{k\choose i+\ell}p^{i}}{{k\choose i-1+\ell}p^{i-1}} =\prod_{1\le i\le x}\frac{(k-i-\ell+1)p}{i+\ell}\le \left(\frac{(k-\ell)p}{\ell+1}\right)^{x}\le \left(\frac{\ell+1}{\ell+1}\right)^{x}=1,$$ using the restriction on $k$ in the first case of \eqref{(kx)p^x-bd}. Then applying \eqref{(kx)p^x-bd} with $\ell=1$ and $x=d-M_n(j)-1$ to \eqref{abs.diff.an} yields
\begin{multline*}
\abs{a_{n,m,(n)}(U_{n},j)-a_{n,m}(U_{n},j)}\\
\le \Cr{absbdJ1}p\left[(M_n(U_n)-d)^+ {\bf 1}\{M_n(U_n)-d\le 2/p+1\}+2^{M_n(U_n)-d}{\bf 1}\{M_n(U_n)-d>2/p+1\}\right]\\
\le \Cr{absbdJ1} p\left[M_n(U_n)+ 2^{M_n(U_n)}{\bf 1}\{M_n(U_n)>2/p\}\right].
\end{multline*}
This same upper bound holds for $j\in J_2$ as well since
\begin{multline*}
\abs{a_{n,m,(n)}(U_{n},j)-a_{n,m}(U_{n},j)}=\abs{-{M_n(U_n)-d\choose d-M_n(j)}p^{d-M_n(j)} (1-p)^{M_n(U_n)+M_n(j)-2d}}\\
\le {M_n(U_n)-d\choose d-M_n(j)}p^{d-M_n(j)}
\end{multline*}
and since $\Cr{absbdJ1}\ge 1$ we obtain
\begin{multline}
E_m\left\{\sum_{1 \le j \le m, j\ne U_n}[a_{n,m,(n)}(U_{n},j)-a_{n,m}(U_{n},j)]\right\}^2 \\
\le E_m\left\{\sum_{j\in J_1\cup J_2}\Cr{absbdJ1} p\left[M_n(U_n)+2^{M_n(U_n)}{\bf 1}\{M_n(U_n)>2/p\}\right]\right\}^2 \\
\le E_m\left\{\Cr{absbdJ1}mp\left[M_n(U_n)+2^{M_n(U_n)}{\bf 1}\{M_n(U_n)>2/p\}\right]\right\}^2 \\
\le 2 (\Cr{absbdJ1}mp)^2\left\{E_mM_n(U_n)^2+ E_m\left[4^{M_n(U_n)}{\bf 1}\{M_n(U_n)>2/p\}\right]\right\}.\label{Ea-diffsq}
\end{multline}
Using that $M_n(U_n)\sim {\cal B}(n,1/m)$  we have
\begin{equation}\label{a:quad-term}
E_mM_n(U_n)^2=(n/m)(1-1/m)+(n/m)^2\le n/m +(n/m)^2 \le 2\left[1+\left(\frac{n}{m}\right)^2\right],\end{equation} using \eqref{xpoly-bd}.
To bound the second expectation in \eqref{Ea-diffsq}, we let $B_{n,s}$ denote a random variable with
distribution ${\cal B}(n,s)$, and note the identity
\begin{equation}\label{tilt-id}
E[w^{B_{n,s}} f(B_{n,s})]=(1-s+sw)^n E f(B_{n,sw/(1-s+sw)}) \qmq{for all $w>0$ and bounded $f$}\end{equation} and the bound
\begin{equation}\label{Hoeff-bd}
P(B_{n,s}>t)\le \exp\left[-2n(t/n-s)^2\right]
\end{equation}
of \cite{Hoeffding63}. Letting $\wtilde{p}=(4/m)/(1+3/m)$ and applying \eqref{tilt-id} with $w=4$ and the bound \eqref{Hoeff-bd}, we have
\begin{multline}
E_m\left[4^{B_{n,1/m}}{\bf 1}\{B_{n,1/m}>2/p\}\right]=(1+3/m)^n P(B_{n,\wtilde{p}}>2/p) \le e^{3n/m} \exp\left[-2n\left(\frac{2}{np}-\wtilde{p}\right)^2\right]\\
= \exp\left[\frac{3n}{m}- \frac{8}{np^2}+\frac{8\wtilde{p}}{p}-2n\wtilde{p}^2 \right]\le \Cl{a:exp-term} \exp\left[\frac{3n}{m}- \frac{8}{np^2}\right]\le \Cr{a:exp-term},\label{a:exp-term}
\end{multline} where in the second-to-last step we used  that $\wtilde{p}/p\le 4$ and $-2n\wtilde{p}^2 \le 0$, and the final step is as follows. The bounds in \eqref{n/m<log<sqrtm} imply that
$$\frac{3n}{m}=3\left(\frac{n}{m}\right)^2\frac{m}{n}\le \frac{m^2}{n},$$
and so
\begin{equation*}
\frac{3n}{m}- \frac{8}{np^2}\le  \frac{m^2}{n}-\frac{8(m-1)^2}{n}= \frac{m^2}{n}-\frac{8m^2}{n}\left(1-\frac{1}{m}\right)^2\le \frac{m^2}{n}-\frac{2m^2}{n}=-\frac{m^2}{n}\le 0,
\end{equation*} implying that the entire term \eqref{a:exp-term} is bounded by $\Cr{a:exp-term}$. Combining this bound with \eqref{a:quad-term} and \eqref{Ea-diffsq} yields
\begin{equation}\label{an:sumjne}
E_m\left\{\sum_{1 \le j \le m, j\ne U_n}[a_{n,m,(n)}(U_{n},j)-a_{n,m}(U_{n},j)]\right\}^2\le \C\left[1+\left(\frac{n}{m}\right)^2\right].
\end{equation}

Now consider the second sum in (\ref{efst_mix_sums}), whose summation index can be further restricted to $i\in J_3\cup J_4$, where
\begin{align*}
J_3&=\{1 \le i \le n, i\ne U_{n}: M_n(i)+M_n(U_n)\ge 2d+1, M_n(U_n)<d\}\\
J_4&=\{1 \le i \le n, i\ne U_{n}: M_n(i)+M_n(U_n)=2d, M_n(U_n)<d\}.
\end{align*} For $i\in J_3$,
\begin{multline}\label{adiffiJ3}
a_{n,m,(n)}(i,U_{n})-a_{n,m}(i,U_{n})={M_n(i)-d\choose d-M_n(U_n)+1}p^{d-M_n(U_n)+1} (1-p)^{M_n(i)+M_n(U_n)-2d-1}\\
-{M_n(i)-d\choose d-M_n(U_n)}p^{d-M_n(U_n)} (1-p)^{M_n(i)+M_n(U_n)-2d}\\
={M_n(i)-d\choose d-M_n(U_n)+1}p^{d-M_n(U_n)} (1-p)^{M_n(i)+M_n(U_n)-2d-1}\left(p-(1-p)\frac{d-M_n(U_n)+1}{M_n(i)+M_n(U_n)-2d}\right)
\end{multline} and, by an argument like \eqref{absbdJ1}, the difference above is bounded in absolute value by
$$p+(1-p)\frac{d-0+1}{(2d+1)-2d}=p+(1-p)(d+1)\le d+1=:\Cl{a-2nd.term}.$$
Applying \eqref{(kx)p^x-bd} with $\ell=2$ and $x=d-M_n(U_n)-1$, we have that \eqref{adiffiJ3} is bounded in absolute value by
\begin{multline*}
\Cr{a-2nd.term}p\left[{M_n(i)-d\choose 2}{\bf 1}\{M_n(i)-d\le 3/p+2\}+2^{M_n(i)-d} {\bf 1}\{M_n(i)-d> 3/p+2\}\right]\\
\le \Cr{a-2nd.term}p\left[M_n(i)^2+2^{M_n(i)} {\bf 1}\{M_n(i)> 3/p\}\right]\le \Cr{a-2nd.term}p\left[M_n(i)^2+2^{M_n(i)} {\bf 1}\{M_n(i)> 2/p\}\right].
\end{multline*}
This same bound holds for $i\in J_4$without the factor $\Cr{a-2nd.term}$, since by \eqref{(kx)p^x-bd} with $\ell=1$ and the same $x$,
\begin{multline*}
\abs{a_{n,m,(n)}(i,U_{n})-a_{n,m}(i,U_{n})}=\abs{-{M_n(i)-d\choose d-M_n(U_n)}p^{d-M_n(U_n)} (1-p)^{M_n(i)+M_n(U_n)-2d}}\\
\le {M_n(i)-d\choose d-M_n(U_n)}p^{d-M_n(U_n)} \le p\left[M_n(i) +2^{M_n(i)} {\bf 1}\{M_n(i)> 2/p\}\right]\\
 \le p\left[M_n(i)^2 +2^{M_n(i)} {\bf 1}\{M_n(i)> 2/p\}\right].\end{multline*}
Hence, as $\Cr{a-2nd.term}\ge 1$,
\begin{multline}\label{Einesq}
E_m\left\{\sum_{i\ne U_n}[a_{n,m,(n)}(i,U_{n})-a_{n,m}(i,U_{n})]\right\}^2
\le E_m\left\{\sum_{i\in J_3\cup J_4}  \Cr{a-2nd.term}p\left[M_n(i)^2+2^{M_n(i)} {\bf 1}\{M_n(i)> 2/p\}\right]\right\}^2\\
\le E_m\left\{\sum_{i=1}^m \Cr{a-2nd.term} p[M_n(i)^2+2^{M_n(i)}{\bf 1}\{M_n(i)>2/p\}]\right\}^2\\
\le (\Cr{a-2nd.term}mp)^2E_m[M_n(1)^2+2^{M_n(1)}{\bf 1}\{M_n(1)>2/p\}]^2\\ \le2(\Cr{a-2nd.term}mp)^2\left\{E_m M_n(1)^4+E_m\left[4^{M_n(1)}{\bf 1}\{M_n(1)>2/p\}\right]\right\}.\end{multline}
Now, combining the bound
\begin{equation*}
E_m M_n(1)^4\le \Cr{Riordin}{_{,4}} [1+(n/m)^4]
\end{equation*}
obtained from \eqref{bin:mom} with $q=4$ and $p=1/m$ with the bound \eqref{a:exp-term}, we have that \eqref{Einesq} can not exceed $\Cl{} [1+(n/m)^4]$. Applying this bound together with \eqref{an:sumjne} in \eqref{efron.stein.n+1} yields
\beas
\var_m(T_{n-1})\le nE_m(T_{n-1}-T_n)^2\le \Cl{vara_n-1}n\left[1+\left(\frac{n}{m}\right)^2+ \left(\frac{n}{m}\right)^4\right]\le 2\Cr{vara_n-1} n\left[1+  \left(\frac{n}{m}\right)^4\right],
\enas
using \eqref{xpoly-bd} in the final inequality. Hence
\begin{equation*}
\var_m\left(\sum_{i\ne j}a_{n,m}(i,j)\right)=\var_m(T_{n})\le 2\Cr{vara_n-1}(n+1)\left[1+\left(\frac{n+1}{m}\right)^4\right]\le \Cr{varsump} n \left[1+\left(\frac{n}{m}\right)^4\right]
\end{equation*} by taking $\Cr{varsump}=2^6\Cr{vara_n-1}$ and using the elementary bound
\begin{multline}\label{elem-bd-nm}
(n+1)\left[1+\left(\frac{n+1}{m}\right)^j\right]=\left(\frac{n+1}{n}\right) n\left[1+\left(\frac{n+1}{n}\right)^j \left(\frac{n}{m}\right)^j\right]\le 2n\left[1+2^j \left(\frac{n}{m}\right)^j\right]\\
\le 2^{j+1}n\left[1+ \left(\frac{n}{m}\right)^j\right]\qmq{for all}j\ge 1,
\end{multline}
thus proving \eqref{varsump}.

To prove \eqref{varsumb} we let $T_n=\sum_{i\ne j} b_{n,m}(i,j)$ and proceed similarly. In view of (\ref{bij.def1}), $b_{n,m}(i,j)=0$ when $M_n(i)+M_n(j) \ge 2d+1$, but since $b_{n,m,(n)}(i,j)$ is calculated when withholding ball $n$, $b_{n,m,(n)}(i,j)$ may be nonzero when $M_n(i)+M_n(j)=2d+1$ and $U_n\in\{i,j\}$, a case we thus allow for in our definition of $K$ below. We have
\begin{align}
T_{n-1}-T_n&=\sum_{j\ne U_{n}}[b_{n,m,(n)}(U_{n},j)-b_{n,m}(U_{n},j)]+\sum_{i\ne U_{n}}[b_{n,m,(n)}(i,U_{n})-b_{n,m}(i,U_{n})]\nonumber\\
&=\sum_{j\in K_1}[b_{n,m,(n)}(U_{n},j)-b_{n,m}(U_{n},j)]+\sum_{i\in K_2}[b_{n,m,(n)}(i,U_{n})-b_{n,m}(i,U_{n})]\label{bn.2sums}
\end{align}
where
\begin{gather*}
 K_1=\{j: (U_n,j)\in K\} \qmq{and} K_2=\{i: (i,U_n)\in K\} \qm{with}\\
K=\{(i,j): i\ne j, M_n(i)+M_n(j)\le 2d+1, M_n(j)>d\}.
\end{gather*}
For any $(i,j) \in K$ we have
\begin{multline}
\frac{{n-d\choose d}{d-M_n(i)\choose M_n(j)-d}}{{n-M_n(i)\choose M_n(j)}} =\frac{{M_n(j)\choose d}{n-M_n(i)-M_n(j)\choose n-2d}}{{n-M_n(i)\choose n-d}}\\
={M_n(j)\choose d}\frac{(n-M_n(i)-M_n(j))_{2d-M_n(i)-M_n(j)}}{(n-M_n(i))_{d-M_n(i)}}\cdot \frac{(d-M_n(i))!}{(2d-M_n(i)-M_n(j))!} \\
\le {2d+1\choose d}\frac{n^{2d-M_n(i)-M_n(j)}}{(n-d+1)^{d-M_n(i)}}\cdot\frac{d!}{1}\le  \Cl{bUj1} \frac{n^{2d-M_n(i)-(d+1)}}{(n-d+1)^{d-M_n(i)}}
=\Cr{bUj1} \frac{n^{d-M_n(i)-1}}{(n-d+1)^{d-M_n(i)}}\\
=\frac{\Cr{bUj1}}{n}\left( \frac{n}{n-d+1}\right)^{d-M_n(i)} \le \frac{\Cr{bUj1}}{n}\left( \frac{2d}{2d-d+1}\right)^{d-M_n(i)} \le \frac{\Cr{bUj1}}{n}\cdot 2^{d-M_n(i)}\le \frac{2^d\Cr{bUj1}}{n},\label{bUj1}
\end{multline}
using that on $K$ we have $M_n(j)\ge d+1$, and that $n \ge 2d$ in the first inequality on the last line.

Now considering $b_{n,m,(n)}(U_n,j)$ for $(U_n,j) \in K$, note that by \eqref{bij.def1}, $n \le 2d$ implies $b_{n,m,(n)}(U_n,j)=0$, as ${n-1-d \choose d}=0$ in this case. Hence we may assume $n-1\ge 2d$. When $(U_n,j) \in K$ when $n$ balls are tossed, then when the $n$th ball $U_n \ne j$ is omitted, $(U_n,j) \in K$ still. Hence in this case (\ref{bUj1}) applies both to $b_{n,m}(U_n,j)$ and $b_{n,m,(n)}(U_n,j)$, yielding
\begin{equation}\label{abs-bK2}
\abs{b_{n,m,(n)}(U_n,j)-b_{n,m}(U_n,j)}\le b_{n,m,(n)}(U_n,j)+b_{n,m}(U_n,j)\le\Cl{bUj2}/n\qmq{for all}j\in K_1.
\end{equation}

For the sum in \eqref{bn.2sums} over $i\in K_2$ we have $b_{n,m}(i,U_n)=0$ if $M_n(i)+M_n(j)=2d+1$, and otherwise the bound (\ref{bUj1}) holds. To consider $b_{n,m,(n)}(i,U_n)$, again assume that $n-1\ge 2d$ since $b_{n,m,(n)}(i,U_n)=0$ otherwise, as before. In addition, if $M_n(U_n)=d+1$ then removing ball $n$ leaves $d$ balls in cell $U_n$, in which case $(i,U_n) \not \in K$, and hence we assume $M_n(U_n)>d+1$. In this case, after removing ball $n$ the pair $(i,U_n)$ remains in $K$, and (\ref{bUj1}) applies. Thus,
\beas 
\abs{b_{n,m,(n)}(i,U_n)-b_{n,m}(i,U_n)}\le b_{n,m,(n)}(i,U_n)+b_{n,m}(i,U_n)\le\Cl{biU2}/n\qmq{for all}i\in K_2,
\enas
and combining this bound with \eqref{abs-bK2} for use in (\ref{efron.stein.n+1}) yields
\begin{multline*}
\var_m(T_{n-1})\le nE_m(T_{n-1}-T_n)^2\\
=nE_m\left\{\sum_{j\in K_1}[b_{n,m,(n)}(U_{n},j)-b_{n,m}(U_{n},j)]+\sum_{i\in K_2}[b_{n,m,(n)}(i,U_{n})-b_{n,m}(i,U_{n})]\right\}^2\\
\le nE_m\left\{\sum_{j\in K_1}\Cr{bUj2}/n+\sum_{i\in K_2}\Cr{biU2}/n\right\}^2\le nE_m\left\{m(\Cr{bUj2}+\Cr{biU2})/n\right\}^2= \Cr{varsumb} \frac{m^2}{n}
\end{multline*} by taking $\Cr{varsumb}=(\Cr{bUj2}+\Cr{biU2})^2$, so
\begin{equation*}
\var_m\left(\sum_{i\ne j}b_{n,m}(i,j)\right)=\var_m(T_{n})\le \Cr{varsumb} \frac{m^2}{n+1}\le \Cr{varsumb} \frac{m^2}{n},
\end{equation*} proving \eqref{varsumb}.

For \eqref{Ass1-3rdsum} we recall expression (\ref{cnm:def}) wherein $q=1-p$, and let $T_n= \sum_{i \ne j} c_{n,m}(i,j)$. Since $c_{n,m,(n)}(i,j)=c_{n,m}(i,j)$ as long as $U_n\not\in \{i,j\}$, we have
$$T_{n-1}-T_n=\sum_{j \ne U_n} [c_{n,m,(n)}(U_n,j)-c_{n,m}(U_n,j)]+\sum_{i \ne U_n} [c_{n,m,(n)}(i,U_n)-c_{n,m}(i,U_n)].$$
Considering the first sum and casing out on whether $M_n(U_n)\le d$ or $M_n(U_n)\ge d+1$,
\begin{multline*}
c_{n,m,(n)}(U_n,j)-c_{n,m}(U_n,j)= {\bf 1}\{M_n(j)=d\}(q^{\abs{M_n(U_n)-d}}-q^{\abs{M_n(U_n)-d-1}})\\
={\bf 1}\{M_n(j)=d\}(pq^{d-M_n(U_n)}{\bf 1}\{M_n(U_n)\le d\}-pq^{M_n(U_n)-d-1}1\{M_n(U_n)\ge d+1\})\\
=p\cdot 1\{M_n(j)=d\}(q^{d-M_n(U_n)}{\bf 1}\{M_n(U_n)\le d\}-q^{M_n(U_n)-d-1}{\bf 1}\{M_n(U_n)\ge d+1\}).
\end{multline*}
Since the term in parentheses is bounded in absolute value by $1$, we have that $$\abs{c_{n,m,(n)}(U_n,j)-c_{n,m}(U_n,j)}\le p\cdot {\bf 1}\{M_n(j)=d\},$$ whence
\begin{multline}
E_m\left(\sum_{j \not = U_n} [c_{n,m,(n)}(U_n,j)-c_{n,m}(U_n,j)]\right)^2\le E_m\left(p\sum_{j=1}^m {\bf 1}\{M_n(j)=d\}\right)^2=p^2E_m(Y_n^2)\\
=p^2(\sigma_{n,m}^2+\mu_{n,m}^2) \le (\C/m)^2\left(\Cr{sigma.over.mu} \mu_{n,m}+ \mu_{n,m}^2\right)\le (\C/m)^2\left(\Cr{sigma.over.mu} m+ m^2\right) \le \Cl{Ass13rdsumA},\label{Ass13rdsumA}
\end{multline} using \eqref{pf:sig<Cmu} and the trivial bound $\mu_{n,m}\le m$. Next,
\begin{multline*}
|c_{n,m,(n)}(i,U_n)-c_{n,m}(i,U_n)| =(1-q^{\abs{M_n(i)-d}})|{\bf 1}\{M_n(U_n)=d+1\}-{\bf 1}\{M_n(U_n)=d\}|\\
            \le 1-q^{\abs{M_n(i)-d}} \le  \abs{M_n(i)-d}\abs{\log q}.
\end{multline*}
Further, since $p=1/(m-1)\le 1/2$ by \eqref{m>=3}, by Taylor series
$$\abs{\log q}=-\log(1-p)\le p+\frac{p^2}{2(1-1/2)^2}=p+2p^2\le \Cl{abslogq}/m,$$ giving
\begin{multline}
E_m\left(\sum_{i \not = U_n} [c_{n,m,(n)}(i,U_n)-c_{n,m}(i,U_n)]\right)^2\le E_m\left(\sum_{i = 1}^m \abs{M_n(i)-d}\abs{\log q}\right)^2\\
\le(\Cr{abslogq}/m)^2m \sum_{i = 1}^mE_m(M_n(i)-d)^2 =\Cr{abslogq}^2 E_m(M_n(1)-d)^2\\
=\Cr{abslogq}^2 [\var_m(M_n(1))+(E_mM_n(1)-d)^2] \le\Cr{abslogq}^2 [\var_m(M_n(1))+2(E_mM_n(1))^2+2d^2]\\
 =\Cr{abslogq}^2 [(n/m)(1-1/m)+2(n/m)^2+2d^2] \le \C{}[1+(n/m)^2]\label{Ass13rdsumB}
\end{multline} by \eqref{xpoly-bd}. Applying \eqref{Ass13rdsumA} and \eqref{Ass13rdsumB} in (\ref{efron.stein.n+1}), we have
$$\var_m(T_{n-1})\le nE_m(T_{n-1}-T_n)^2\le \Cl{varsumc:n-1} n\left[1+\left(\frac{n}{m}\right)^2\right],$$ so
\begin{equation*}
\var_m\left(\sum_{i\ne j}c_{n,m}(i,j)\right)=\var_m(T_n)\le \Cr{varsumc:n-1} (n+1)\left[1+\left(\frac{n+1}{m}\right)^2\right]\le \Cr{Ass1-3rdsum}n\left[1+\left(\frac{n}{m}\right)^2\right]\end{equation*}
by taking $\Cr{Ass1-3rdsum}=8\Cr{varsumc:n-1}$ and using \eqref{elem-bd-nm}. This proves \eqref{Ass1-3rdsum} and thus concludes the proof of the lemma.\bbox

\section{Discussion}\label{discussion}
\begin{table}[htp]
\caption{Asymptotic domains, wherein $\tau_d$ and $\varphi_d$ are given by \eqref{tau} and $f(n,m) \sim g(n,m)$ denotes  $f(n,m)=(1+o(1))g(n,m)$. That $\lim \sigma_{n,m}^2/\mu_{n,m}$  is strictly positive in the central domain follows from Lemma~\ref{r1.large}, Part~\ref{infvarphi>0}.}
\begin{center}
\begin{tabular}{c|c|c|c}
Definition&$\mu_{n,m}$&$\sigma_{n,m}^2/\mu_{n,m}$&Asymptotic  Distribution\\
&&&of $Y_n$ under $P_m$\\\hline\hline
\multicolumn{4}{c}{Left-hand domain}\\
$n/m\To 0$&$\To\mu$&$\To 1$&$\mbox{Poi}(\mu)$\\
$\mu_{n,m}\To\mu\in(0,\infty)$&&&\\\hline
\multicolumn{4}{c}{Left-hand intermediate domain}\\
$n/m\To 0$&$\sim m\tau_d(n/m)$&$\To 1$&$N(\mu_{n,m},\sigma_{n,m}^2)$\\
$\mu_{n,m}\To\infty$&$\To\infty$&&\\\hline
\multicolumn{4}{c}{Central domain}\\
$n/m\To\rho\in(0,\infty)$&$\sim m\tau_d(\rho)$&$\To \varphi_d(\rho)\in(0,1)$&$N(\mu_{n,m},\sigma_{n,m}^2)$\\\hline
\multicolumn{4}{c}{Right-hand intermediate domain}\\
$n/m\To\infty$&$\sim m\tau_d(n/m)$&$\To 1$&$N(\mu_{n,m},\sigma_{n,m}^2)$\\
$\mu_{n,m}\To \infty$&$\To \infty$&&\\\hline
\multicolumn{4}{c}{Right-hand domain}\\
$n/m\To\infty$&$\To \mu$&$\To 1$&$\mbox{Poi}(\mu)$\\
$\mu_{n,m}\To\mu\in(0,\infty)$&&&\\\hline
\end{tabular}
\end{center}
\label{table:domains}
\end{table}

Theorem \ref{thm:main} is applied in \citet{Go12} to obtain bounds on the normal approximation for the number of vertices in the Erd\H{o}s-R\'enyi random graph of a given degree. Although the graph degree and occupancy problems have some features in common, they also differ in a number of significant ways. On balance, the occupancy problem is the more difficult of the two for the following reasons.

First, the term $\Psi_{n,m}^2$ required by Condition \ref{Psinisrootn} of Theorem \ref{thm:main}, the variance of the conditional expectation of the difference $Y_n^s-Y_n$, is harder to compute for the occupancy problem. In particular, in the graph degree problem one can make a direct bound on this term, but here we appear to be forced to rely instead on the use of the Efron-Stein inequality.

Another significant difference between these two problems is that for graph degree the removal of a vertex leaves the connectivity of the remaining graph unaffected, while the parameters of the occupancy problem that results after the removal of an urn depends on the number of balls that urn contained. In particular, even if the removed vertex in the graph degree problem was connected to all other vertices the reduced graph remains non-trivial, in contrast to the `parallel' situation of removing an urn which contains all balls in the occupancy problem. As a result, though the graph degree problem is indexed by the number of vertices, and the variable of interest is a count over those same vertices, here we index by the number of balls, while the count is a sum over urns. The choice is driven by the fact that Condition \ref{assumption:condlawYnvgivenK} is concerned only with reduced problems of sizes $n-L_n$ that satisfy $L_n \le s_n$. And in the occupancy problem, limiting the number of urns that are removed when forming the subproblem does not guarantee that the reduced problem will be non-trivial, but limiting the number of balls removed does.

A third important difference is that in the graph degree problem, we consider a graph with $n$ vertices and connectivity $\theta/(n-1)$, and the reduced problem is on the graph with one vertex removed. There, choosing the parameter space to be $\Theta_n=(0,b] \cap (0,n-1)$ for some large $b$ yields that $\theta \in \Theta_n$ implies $\psi_{n,\theta} \in \Theta_{n-1}$ where $\psi_{n,\theta}=(n-2)\theta/(n-1)$, as required by Condition \ref{assumption:condlawYnvgivenK} of Theorem \ref{thm:main}. That each parameter space $\Theta_n$ is a subset of the same bounded interval $(0,b]$ simplifies a number of the computations and bounds. For the occupancy problem we have taken $\Theta_n$ to be unbounded for the following reason. When an empty cell is removed to form the reduced problem, the differences $n/(m-1)-n/m$ of the ratio of balls to urns equals $n/(m(m-1))$. In the central domain this ratio behaves like $1/m$, summing to the divergent harmonic series. On the other hand, though we appear forced to deal with the case where $\Theta$ is unbounded, here we obtain results in asymptotic domains in addition to the central one.

Although we state our main result, Theorem~\ref{thm:occ}, as a uniform bound holding for all $n\ge d$ and $m\ge 2$, the occupancy problem has been classically studied asymptotically as $n, m\To\infty$, such as by \citet{Kolchin78}, who define the five asymptotic domains given in Table~\ref{table:domains}, which together give an essentially complete asymptotic picture of the $n,m\To\infty$ asymptotic with $n/m$ varying from $0$ to $\infty$.
\citet[][Theorem~4, p.\ 68]{Kolchin78} also show that, in the uniform occupancy problem, $Y_n$ has limiting normal distribution in exactly the domains in which $\sigma_{n,m}\To\infty$, i.e., in the left-hand intermediate, central, and right-hand intermediate domains.  Except for a small portion of the latter, our Berry-Esseen type bound in Theorem~\ref{thm:occ} provides convergence to the normal in these domains as well: The left-hand intermediate and central domains are covered by Corollary~\ref{cor:occ}, and the right-hand intermediate domain is addressed in the following.

\begin{corollary}\label{cor:RHIdomain}
Let $$\delta_{n,m}=\log m+d\log\log m-n/m.$$ If $n,m\To\infty$ in such a way that $n/m\To\infty$, $\mu_{n,m}\To\infty$, and
\begin{equation}\label{supn/mlogm<1}
\lim_{n,m\To\infty} \left(\frac{\delta_{n,m}}{\log\log m}\right)>6,
\end{equation} then $r_{n,m}\To\infty$ and, in particular, $$\sup_{z \in \mathbb{R}}\left|P\left( W_{n,m} \le z \right) -P(Z \le z)\right|\To 0.$$
\end{corollary}

\proof If \eqref{supn/mlogm<1} holds then there is $\eps>0$ such that, for all $n, m$ sufficiently large,
\begin{equation}\label{n/m.upper}
\frac{\delta_{n,m}}{\log\log m}\ge 6+\eps,\qmq{or equivalently}n/m\le \log m+(d-6-\eps)\log\log m.\end{equation} We will use below that $\log[x(1+o(1))]=\log x+o(1)$. Using Table \ref{table:domains},
\begin{multline}\label{r^2-right}
\log r_{n,m}^2=\log\left\{\frac{\sigma_{n,m}^2}{[1+(n/m)^3]^2}\right\}=\log\left\{(1+o(1))\frac{m\tau_d(n/m)}{(n/m)^6}\right\}\\
=\log\left\{\frac{m(n/m)^{d}e^{-n/m}}{(n/m)^6 d!}\right\}+o(1)=\log m+(d-6)\log(n/m)-n/m-\log(d!)+o(1)\\
=\log m+(d-6)\log(n/m)-n/m+O(1).
\end{multline}
Noting that $x\mapsto (d-6)\log x-x$ is decreasing for $x>(d-6)^+$, for $n,m$ large enough so that $n/m\ge (d-6)^+$ and \eqref{n/m.upper} holds, by \eqref{r^2-right} we have
\begin{multline*}
\log r_{n,m}^2\ge \log m+(d-6)\log\left[\log m+(d-6-\eps)\log\log m\right]\\
-\left[\log m+(d-6-\eps)\log\log m\right]+O(1)\\
=(d-6)\log[(1+o(1))\log m]- (d-6-\eps)\log\log m +O(1)\\
=(d-6)\log\log m+o(1)- (d-6-\eps)\log\log m +O(1) =\eps \log\log m +O(1)\To\infty.
\end{multline*}\bbox

An example of a regime satisfying the hypothesis of Corollary~\ref{cor:RHIdomain} is
\bea \label{non.vac}
n=\lfloor m \left( \log m+(d-a)\log\log m \right) \rfloor,\quad a>6.
\ena
Then $\log\mu_{n,m}=a\log\log m+O(1)\To\infty$ by (\ref{occ:mv}) and (\ref{non.vac}), and $\delta_{n,m}/(\log\log m)\To a$, so \eqref{supn/mlogm<1} is satisfied.

Although \eqref{supn/mlogm<1} does not cover all of the right-hand intermediate domain, the missing part is small since it follows from $n/m\To\infty$ and $\mu_{n,m}\To\infty$ that $\delta_{n,m} \To\infty$ (see, e.g., \citet[][p.~41]{Kolchin78}). Since \eqref{supn/mlogm<1} implies that $\delta_{n,m}=a(\log\log m)(1+o(1))$ for some $a>6$, the only portion of the right-hand intermediate domain  in which $r_{n,m}\not\To\infty$ but $Y_n$ still converges to normal is the narrow asymptotic where $\delta_{n,m} \To\infty$ but $\delta_{n,m}\le 6(\log\log m)(1+o(1))$. We note, however, that the limiting factor in $r_{n,m}$ that determines our bound arises when bounding $\Psi_n^2$, a term that also appears when computing a bound to the normal in the weaker Wasserstein metric using size bias couplings, such as the bound obtained by applying Corollary 2.2 and Construction 3A of \citet{ChRo10}. Hence this behavior appears to be unrelated to any aspect of our method that pertains to bounding the Kolmogorov distance, and the inductive method in particular.

Lastly we remark that extensions of the present work to the case where the cell probabilities are non-uniform is of additional interest, and may likely also be approached with the use of Theorem \ref{thm:main}.

\section*{Appendix}

\subsection*{Proof of Lemma \ref{r1.large}}
Part~\ref{Theta.finite}: By (\ref{r.min.subset}), it suffices to show (\ref{sig-to-0-m}), as then for any $r_1>0$ the set on the right hand side of (\ref{r.min.subset}) is finite. By (\ref{occ:mv}) and $d \ge 2$ we have
\beas
\limsup_{m \rightarrow \infty} \mu_{n,m} = \limsup_{m \rightarrow \infty} {n \choose d}\frac{1}{m^{d-1}}\left(1-\frac{1}{m}\right)^{n-d} \le {n \choose d} \limsup_{m \rightarrow \infty} \frac{1}{m^{d-1}}=0,
\enas
and similarly
\beas
\limsup_{m \rightarrow \infty} m(m-1){n \choose d,d,n-2d}\frac{1}{m^{2d}}\left( 1 - \frac{2}{m} \right)^{n-2d} \le
{n \choose d,d,n-2d} \limsup_{m \rightarrow \infty} \frac{1}{m^{2(d-1)}}=0.
\enas
Hence (\ref{sig-to-0-m}) holds by (\ref{varYn}).

Part~\ref{sig2.le.mu}: As the mean $\mu_{n,m}$ is positive over the range of $n$ and $m$ considered, we equivalently show that the ratio
\beas
\frac{\sigma_{n,m}^2}{\mu_{n,m}} = 1-\mu_{n,m} + \frac{{n \choose d,d,n-2d}\left( 1 - \frac{2}{m} \right)^{n-2d}}{{n \choose d}\left(1-\frac{1}{m}\right)^{n-d}} \left( \frac{m-1}{m^d} \right)
\enas
is bounded. For $d \le n < 2d$ or $m=2$ the result is clear, as
\beas
\frac{\sigma_{n,m}^2}{\mu_{n,m}} = 1-\mu_{n,m}  \le 1.
\enas
For $n \ge 2d$ and $m \ge 3$, we obtain
\bea
\lefteqn{\frac{\sigma_{n,m}^2}{\mu_{n,m}}-1=\frac{{n \choose d,d,n-2d}\left( 1 - \frac{2}{m} \right)^{n-2d}}{{n \choose d}\left(1-\frac{1}{m}\right)^{n-d}} \left( \frac{m-1}{m^d} \right) - {n \choose d}\frac{1}{m^{d-1}}\left(1-\frac{1}{m}\right)^{n-d}}\nn \\
&&= \frac{(n)_{2d}\left( 1 - \frac{2}{m} \right)^{-d}}{d!(n)_d} \left( \frac{m-1}{m^d} \right)
\left( 1-\frac{1}{m-1} \right)^{n-d}
- \frac{(n)_d}{d!}\frac{1}{m^{d-1}}\left(1-\frac{1}{m}\right)^{n-d}\nn \\
&&=\frac{n^d}{d! m^{d-1}} \left(\left( 1 - \frac{2}{m} \right)^{-d}\left( 1-\frac{1}{m-1} \right)^{n-d}- \left(1-\frac{1}{m}\right)^{n-d} \right)\label{two.taylor}\\
&&-\frac{(n)_{2d}\left( 1 - \frac{2}{m} \right)^{-d}}{d!(n)_dm^d}
\left( 1-\frac{1}{m-1} \right)^{n-d}\nn \\
&&+O\left(\left(\frac{n}{m}\right)^{d-1}\right)\left(\left( 1 - \frac{2}{m} \right)^{-d}\left( 1-\frac{1}{m-1} \right)^{n-d} + \left(1-\frac{1}{m}\right)^{n-d} \right). \nn
\ena
For the first term in (\ref{two.taylor}), expanding $(1-x)^{-d}$ around zero and using $d \ge 3$ yields
\beas
\left(1-\frac{2}{m} \right)^{-d} = 1 + \frac{2d}{m} + R_{1,m} \qmq{with} |R_{1,m}| \le \frac{2d(d+1)3^{d+2}}{m^2},
\enas
and similarly expanding $(1-x)^{n-d}$ around $1/m$ yields
\beas
\left(1-\frac{1}{m-1}\right)^{n-d} =
\left(1-\frac{1}{m}\right)^{n-d} + R_{2,m}
\enas
with
\beas
|R_{2,m}| \le \frac{n-d}{m(m-1)}\left(1-\frac{1}{m}\right)^{n-d-1} \le \frac{n}{m(m-1)}e^{-(n-d-1)/m}.
\enas
Hence, we may bound the first term in (\ref{two.taylor}) by $n^d/d!m^{d-1}$ times
\beas
\lefteqn{\left| \left(1+\frac{2d}{m}+R_{1,m}\right) \left(\left(1-\frac{1}{m}\right)^{n-d}+R_{2,m} \right) -
\left(1 -\frac{1}{m} \right)^{n-d} \right|}\\
&&=\left| \left( 1 + \frac{2d}{m} \right)R_{2,m}+\frac{2d}{m}\left( 1- \frac{1}{m} \right)^{n-d}
+R_{1,m}\left( \left(1-\frac{1}{m}\right)^{n-d}+R_{2,m}\right) \right|\\
&&\le \frac{e^{-n/m}}{m}\left( \left( 1+\frac{2d}{m} \right)\frac{ne^{(d+1)/m}}{m-1}+ 2d e^{d/m}+ \frac{2d(d+1)3^{d+2}}{m} \left(e^{d/m}+\frac{ne^{(d+1)/m}}{m(m-1)} \right)\right)\\
&&\le \frac{\Cl{dep.on.d} e^{-n/m}}{m}\left(1+ \frac{n}{m}\right).
\enas
where we have used bounds such as $e^{d/m} \le e^{d/3}$ for $m \ge 3$, and where $\Cr{dep.on.d}$ is a constant depending only on $d$. Hence the first term in (\ref{two.taylor}) can be no greater than
\beas
\frac{\Cr{dep.on.d}}{d!} e^{-n/m}\left( \frac{n}{m}\right)^d \left( 1+ \frac{n}{m} \right) \le \frac{\Cr{dep.on.d}}{d!} \sup_{x \ge 0}e^{-x} x^d (1+x)
\enas
for all $n \ge 2d,m \ge 3$.

The next term in (\ref{two.taylor}) is also bounded, as
\beas
\left| \frac{(n)_{2d}\left( 1 - \frac{2}{m} \right)^{-d}}{d!(n)_dm^d}
\left( 1-\frac{1}{m-1} \right)^{n-d} \right| \le \C \left( \frac{n}{m} \right)^d e^{-n/m},
\enas
as are the final terms, in view of
\beas
\left(\frac{n}{m}\right)^{d-1}\left(\left( 1 - \frac{2}{m} \right)^{-d}\left( 1-\frac{1}{m-1} \right)^{n-d} + \left(1-\frac{1}{m}\right)^{n-d} \right) \le \C \left(\frac{n}{m}\right)^{d-1} e^{-n/m}.
\enas

Part~\ref{mu,sig<n}: Using \eqref{muupbd},
\beas
\frac{\mu_{n,m}}{n}\le \frac{\tau_d(n/m)e^{d/m}}{n/m}=\frac{e^{-n/m}(n/m)^{d-1}e^{d/m}}{d!}\le \frac{e^d}{d!} \sup_{x>0} e^{-x}x^{d-1}
\enas so taking this to be $\Cr{mu<n}$ suffices to show the first claim.  The second now follows from
Part~\ref{sig2.le.mu}.

Part~\ref{infvarphi>0}: By differentiating and factoring we have
\begin{align*}
d!\cdot \varphi_d'(x)&=x^{d-2}e^{-x}[x^3-3dx^2+d(3d-1)x-d^2(d-1)] \\
&=x^{d-2}e^{-x}[x-(d-\sqrt{d})][x-d][x-(d+\sqrt{d})],
\end{align*} and by considering the sign of this derivative we see that
$\inf_{x>0}\varphi_d(x)=\min \varphi_d(d\pm\sqrt{d})$, which we now show is positive. Letting $y$ denote either positive value $d\pm\sqrt{d}$, note that
\begin{align}
\frac{y^{d+1}}{d!}=\frac{y^{d-1}(d-1-y)^2}{(d-1)!}\cdot \frac{y^2}{d(d-1-y)^2} &=\frac{y^{d-1}(d-1-y)^2}{(d-1)!}\cdot \frac{(d\pm \sqrt{d})^2}{d(-1\mp\sqrt{d})^2}\nonumber \\
&=\frac{y^{d-1}(d-1-y)^2}{(d-1)!}\cdot 1.\label{dpmsqrtd}
\end{align} Noting also that $\sum_{d'=0}^\infty e^{-y}y^{d'}(d'-y)^2/d'!=y$ by considering the variance of a Poisson random variable with mean $y$, we have
\begin{align*}
\varphi_d(y)&=1-\tau_d(y)-\tau_d(y)(y-d)^2/y \\
&=1-\frac{1}{y}\left(\frac{e^{-y}y^{d+1}}{d!}+\frac{e^{-y}y^{d}(y-d)^2}{d!}\right)\\
&=1-\frac{1}{y}\left(\frac{e^{-y}y^{d-1}(d-1-y)^2}{(d-1)!}+\frac{e^{-y}y^{d}(d-y)^2}{d!}\right)\qm{(by \eqref{dpmsqrtd})}\\
&>1-\frac{1}{y}\sum_{d'=0}^\infty \frac{e^{-y}y^{d'}(d'-y)^2}{d'!}\\
&=0.
\end{align*}

Part~\ref{nmToinf}: If the claim fails then there are $n^*, m^*<\infty$ and sequences $r_j\To\infty$ and $(n_j,m_j)$ such that, for each $j$, $\sigma_{n_j,m_j}^2\ge r_j$ but $n_j\le n^*$ or $m_j\le m^*$. As at least one of the previous two inequalities must hold for infinitely many $j$, by considering subsequences we may assume that
\begin{equation}
n_j\le n^*\;\mbox{for all $j$,}\qmq{or} m_j\le m^*\;\mbox{for all $j$.}\label{n-not-inf}
\end{equation}
If the first case of \eqref{n-not-inf} holds we have
\begin{equation}\label{maxs.to.inf}
\infty=\lim_{j \rightarrow \infty} r_j=\lim_{j \rightarrow \infty}\sigma_{n_j,m_j}^2=\liminf_{j \rightarrow \infty}\sigma_{n_j,m_j}^2\le\liminf_{j \rightarrow \infty}\left(\max_{n\le n^*}\sigma_{n,m_j}^2\right).
\end{equation} If it were that $\overline{m}:=\sup_j m_j=\infty$, then there would be a subsequence $j_k\To\infty$ on which $m_{j_k}\To\infty$, hence $\sigma_{n,m_{j_k}}^2\To 0$ for each fixed $n$ by Part~\ref{Theta.finite}, and by taking the maximum over a finite set,
\begin{equation}\label{0=lim>liminf}
0=\lim_{k\To\infty}\max_{n\le n^*}\sigma_{n,m_{j_k}}^2 \ge \liminf_{j \rightarrow \infty}\left(\max_{n\le n^*}\sigma_{n,m_j}^2\right),
\end{equation} this inequality because the $\liminf$ of a sequence is always at least as small as the limit of a subsequence. \eqref{0=lim>liminf} would be a contradiction of  (\ref{maxs.to.inf}), leaving only the possibility that $\overline{m}<\infty$, which also leads to a contradiction since, using the second equality in \eqref{maxs.to.inf},
$$\infty= \lim_j \sigma_{n_j,m_j}^2\le \sup_j\sigma_{n_j,m_j}^2\le \sup_{n\le n^*,\; m\le \overline{m}}\sigma_{n,m}^2<\infty,$$ as the last supremum is taken over a finite set. If the second case of \eqref{n-not-inf} holds then for $j$ large enough so that $r_j>(m^*)^2$, since $Y_n\le m$ under $P_m$ we have
$$(m^*)^2<r_j\le \sigma_{n_j,m_j}^2=\var_{m_j}(Y_{n_j})\le E_{m_j}(Y_{n_j}^2)\le m_j^2\le (m^*)^2,$$ again a contradiction.

Part~\ref{n/mlogm}: If the claim fails then there is $\eps>0$ and sequences $r_j\To\infty$ and $(n_j,m_j)$ such that, for each $j$, $\sigma_{n_j,m_j}^2\ge r_j$ but $n_j/m_j> (1+\eps)\log m_j$. By Part~\ref{nmToinf}, taking subsequences if necessary, we can assume that $n_j, m_j\To\infty$, and that for all $j$
\begin{equation}\label{n/mlogmd}
n_j/m_j> (1+\eps)\log m_j>d.\end{equation} Then, using Part~\ref{sig2.le.mu} we obtain that for all $j$ sufficiently large
\begin{multline*}
\frac{r_j}{\Cr{sigma.over.mu}}\le \frac{\sigma_{n_j,m_j}^2}{\Cr{sigma.over.mu}} \le \mu_{n_j,m_j} =m_j{n_j \choose d}\frac{1}{m_j^d}\left(1-\frac{1}{m_j}\right)^{n_j-d}\\
\le m_j\left(\frac{n_j}{m_j}\right)^d \exp[-(n_j-d)/m_j]\le m_j\left(\frac{n_j}{m_j}\right)^d e^{-n_j/m_j} e^{d}.
\end{multline*} Using \eqref{n/mlogmd} and that the function $x\mapsto x^d e^{-x}$ is decreasing for $x>d$, we have
$$\frac{r_j}{\Cr{sigma.over.mu}}\le m_j\left(\frac{n_j}{m_j}\right)^d e^{-n_j/m_j} e^{d}\le m_j\left((1+\eps)\log m_j\right)^d m_j^{-(1+\eps)} e^{d}\le m_j^{-\eps}\left((1+\eps)\log m_j\right)^d e^{d},$$ giving the contradiction
$$\infty=\lim_{j \rightarrow \infty} r_j\le \Cr{sigma.over.mu}(1+\eps)^de^d\lim_{j \rightarrow \infty} \left[m_j^{-\eps}\left(\log m_j\right)^d\right]=0.$$

Part~\ref{sig/mu>C}: Set $\eps=1/4$ and, using \eqref{Kolsig} and Parts~\ref{nmToinf} and \ref{n/mlogm}, choose $r_1>0$ large enough so that any $n,m$ satisfying $\sigma_{n,m}^2\ge r_1$ also satisfy
$$\sigma_{n,m}^2\ge m\tau_d(n/m)\varphi_d(n/m)(1-\eps)\qmq{and}e^{-d/m}\ge 1-\eps.$$ Then, using \eqref{muupbd},
$$\frac{\sigma_{n,m}^2}{\mu_{n,m}}\ge\frac{m\tau_d(n/m)\varphi_d(n/m)(1-\eps)}{m\tau_d(n/m)e^{d/m}} \ge \varphi_d(n/m)(1-2\eps)\ge (1/2)\inf_{x>0}\varphi_d(x),$$ so taking $\Cr{sig/mu>C}$ to be this last suffices.
\bbox

\def\cprime{$'$}







%

%
\end{document}